\documentclass[manuscript,nonacm]{acmart}

\usepackage{graphicx}
\usepackage{amsmath}
\usepackage{amsthm}
\usepackage{booktabs}
\usepackage{tabularx}
\usepackage{algpseudocode}
\usepackage{physics}
\usepackage{backnaur}
\usepackage{subcaption}
\usepackage{hyperref}
\usepackage{tikz}

\usepackage{algorithmicx}
\usepackage{algorithm}
\usepackage{enumitem}
\usepackage{booktabs}
\usepackage{xcolor}
\colorlet{lightred}{red!50}
\colorlet{lightblue}{blue!50}
\begin{document}
\title{Evolving Generalizable Multigrid-Based Helmholtz Preconditioners with Grammar-Guided Genetic Programming}
\author{Jonas Schmitt}
\email{jonas.schmitt@fau.de}
\orcid{0000-0002-8891-0046}
\affiliation{%
	\institution{Friedrich-Alexander-Universität Erlangen-Nürnberg}
	\streetaddress{Cauerstraße 11}
	\city{Erlangen}
	\country{Germany}
	\postcode{91058}
}

\author{Harald Köstler}
\email{harald.koestler@fau.de}
\orcid{0000-0002-6992-2690}
\affiliation{%
	\institution{Friedrich-Alexander-Universität Erlangen-Nürnberg}
	\streetaddress{Cauerstraße 11}
	\city{Erlangen}
	\country{Germany}
	\postcode{91058}
}

\begin{abstract}
Solving the indefinite Helmholtz equation is not only crucial for the understanding of many physical phenomena but also represents an outstandingly-difficult benchmark problem for the successful application of numerical methods.
Here we introduce a new approach for evolving efficient preconditioned iterative solvers for Helmholtz problems with multi-objective grammar-guided genetic programming.
Our approach is based on a novel context-free grammar, which enables the construction of multigrid preconditioners that employ a tailored sequence of operations on each discretization level.
To find solvers that generalize well over the given domain, we propose a custom method of successive problem difficulty adaption, in which we evaluate a preconditioner's efficiency on increasingly ill-conditioned problem instances.
We demonstrate our approach's effectiveness by evolving multigrid-based preconditioners for a two-dimensional indefinite Helmholtz problem that outperform several human-designed methods for different wavenumbers up to systems of linear equations with more than a million unknowns.
\end{abstract}

\maketitle

\section{Introduction}\label{sec:intro}
Automated algorithm design is a long-standing challenge in artificial intelligence (AI) and has two essential goals: \emph{Generalization} and \emph{efficiency}.
Thus, the designed algorithm should not only produce the correct output for arbitrary inputs, but the goal is also to achieve better performance than methods designed by a human expert.
In this work, we aim to demonstrate that this goal is attainable for the indefinite Helmholtz equation, an important benchmark problem from the domain of numerical mathematics.
The Helmholtz equation frequently arises in the study of physical phenomena, such as electromagnetics~\cite{salon1999numerical} and acoustics~\cite{pierce2019acoustics}, and is given by the linear partial differential equation (PDE)
\begin{equation}
	-\nabla ^{2}u - k^{2}u = f,
	\label{eq:helmholtz}
\end{equation}
where $\nabla^2$ is the Laplace operator, $k$ the \emph{wavenumber}, and $f$ the source term. 
In general, the analytic solution $u$ of this equation is unknown, which necessitates the use of numerical methods.
Unfortunately, for large wavenumbers, the system of linear equations that arises from the resulting discretization becomes indefinite and highly ill-conditioned, which means that even small perturbations, for instance due to numerical inaccuracies, have a dramatic effect on the overall error of the computed approximation.
As a consequence, the efficient solution of the indefinite Helmholtz equation is still an open challenge in numerical mathematics~\cite{ernst2012difficult,erlangga2008advances,antoine2016integral}.
Even though various methods for solving this equation have been proposed, many of them fail to generalize over different ranges of wavenumbers and are, thus, only limited to certain problem instances.
Therefore, the design of an efficient Helmholtz solver is not only of great significance for many real-world problems~\cite{versteeg1994marmousi,martin2006marmousi2,billette20052004} but also represents a challenging benchmark for the application of AI-based methods.

Data-driven~\cite{li2020fourier,kovachki2021neural} and physics-informed~\cite{karniadakis2021physics,raissi2019physics} machine learning models have recently achieved significant progress in solving PDEs.
While these approaches have shown competitive or even improved performance compared to classical solvers, their behavior on unseen problems is often difficult to predict, and hence generalization is only possible to a limited degree~\cite{li2021physics}.
Additionally, many of these models require an enormous amount of training data, whose generation still depends on the utilization of conventional numerical methods.
An alternative approach is the application of AI-based optimization methods to improve the efficiency of an existing solver.
In contrast to a trained machine learning model, numerical methods can be formulated in the language of mathematics in a problem-independent manner, which significantly facilitates their generalizability.
Furthermore, as this formulation can be understood by a human expert, it is possible to profoundly analyze their behavior based on existing domain knowledge and expertise.
A numerical method that has achieved considerable success in solving Helmholtz problems~\cite{oseikuffuor2010preconditioning,erlangga2008advances,gander2019solvers}, as well as many other complex PDEs~\cite{benzi2005numerical}, is the acceleration of a slowly converging iterative method with the application of a so-called preconditioner $M$~\cite{benzi2002preconditioning}.
This approach is based on the idea of considering the modified system of linear equations 
\begin{equation}
	A M^{-1} \hat{u} = f,
	\label{eq:helmholtz-preconditioned}
\end{equation}
where $A$ represents the discretized operator of the original system.
An idea first proposed in \cite{erlangga2004preconditioner} for preconditioning Helmholtz systems is the choice of the $M$ as a complex-shifted version of $A$, which results in
\begin{equation}
	M = -\nabla ^{2} - (k^{2} + \varepsilon i).
\end{equation}
In \cite{cocquet2017shift} it is shown that the choice of $\varepsilon$ represents a compromise between the resulting system's solvability and the effectiveness of the preconditioner.
After choosing $M$, the main requirement is to have an efficient method for solving the system 
\begin{equation}
	M u = \hat{u},
	\label{eq:preconditioning-system}
\end{equation}
such that $u \approx M^{-1} \hat{u}$.
Multigrid methods are a class of numerical methods for solving discretized PDEs that are ideally suited for this task~\cite{trottenberg2000multigrid,hackbusch2013multi,briggs2000multigrid}.
If properly constructed, these methods achieve $h$-independent convergence while only requiring $\mathcal{O}(n)$ operations, which means that the number of iterations required for solving a system of linear equations with $n$ unknowns is independent of the discretization width $h$.
While a suitable preconditioning matrix $M$ can often be obtained by analyzing the properties of the system matrix $A$, constructing an efficient multigrid method for its inversion is usually less intuitive and can, thus, be considered a problem of optimal algorithm design.
Since their invention by Federenko and Brandt~\cite{fedorenko1962relaxation,brandt1977multi}, multigrid solvers have been designed predominantly by hand.
Only in recent decades, the automated optimization of these methods has become an active field of research.
\paragraph{Related Work on Multigrid Solver Design}
In principle, a multigrid method is characterized by a finite number of components and design choices that determine its computational structure:
The choice of the \emph{smoother}, \emph{prolongation} and \emph{restriction} operator, \emph{coarse grid solver} and \emph{cycle type}~\cite{trottenberg2000multigrid,briggs2000multigrid}.
A common approach to automate multigrid solver design is to formulate the task as a discrete optimization problem, which is then solved, for instance, using an evolutionary algorithm~\cite{oosterlee2003genetic}, branch-and-bound~\cite{thekale2010optimizing}, or minimax approach~\cite{brown2021tuning}.
A different direction is the application of machine learning methods either to optimize the individual components of a multigrid solver, as the prolongation operator~\cite{greenfeld2019learning,katrutsa2020black,luz2020learning} and smoother~\cite{huang2021learning}, or by replacing certain steps within the method altogether by a machine learning system~\cite{taghibakhshi2021optimization}.
All these approaches have in common that they consider a multigrid method's algorithmic structure immutable. Each step of the method employs a fixed sequence of operations in the form of a particular cycle.
Multigrid cycles are commonly classified into three different categories, V-, W-, and F-cycles, where each cycle type exhibits a distinct computational pattern that represents a compromise between the amount of work performed and the expected speed of convergence~\cite{trottenberg2000multigrid}.
In~\cite{schmitt2021evostencils,schmitt2020constructing} we have proposed a context-free grammar that allows alternating each step of a multigrid solver independently.
Consequently, the search space produced by this grammar includes methods that do not fit into any of the known categories.
While until recently, solvers of such unconventional structure have not been considered, in~\cite{schmitt2021evostencils} it could be demonstrated that multigrid methods evolved by a grammar-based genetic programming approach can achieve higher efficiency in solving certain PDEs than traditional variants.
However, in contrast to the indefinite Helmholtz equation, the PDEs considered in this work can already be efficiently solved with standard multigrid cycles without requiring any further optimization.
Furthermore, while we have demonstrated that the solvers obtained by this approach can also function for larger problem instances than those considered within the search, a systematic approach to generalize a multigrid method to a family of problem instances that share common characteristics is still missing.
To overcome these limitations and extend the context-free grammar introduced in~\cite{schmitt2021evostencils,schmitt2020constructing} to the domain of multigrid preconditioners, we make the following contributions. 
\paragraph{Our Contributions}
We introduce a multi-objective grammar-guided evolutionary search method for finding multigrid preconditioners that generalize well over a sequence of increasingly difficult problem instances and demonstrate its effectiveness by evolving preconditioners for the discretized Helmholtz equation with increasingly high wavenumbers.
\begin{itemize}
	\item To apply our evolutionary search method to the domain of multigrid preconditioners, we adapt the class of context-free grammars presented in~\cite{schmitt2021evostencils,schmitt2020constructing} such that the generated methods can be integrated into an existing iterative solver as a preconditioner. 
	To our knowledge, this is the first formal system that enables the application of grammar-guided genetic programming (GGGP)~\cite{mckay2010grammar,whigham1995grammatically} to the design of preconditioned iterative solvers in a generalizable way.
	\item Since our grammar-based representation of multigrid preconditioners is problem-size independent, each method can be ported and applied to similar problem instances without the need to adapt its internal structure. 
	\item Our evolutionary search method is based on classical tree-based GGGP but copes with the high computational demands for solving PDEs numerically by combining multi-objective optimization with a custom method of successive problem difficulty adaption based on the $h$-independent convergence of multigrid methods.
	\item We demonstrate that our implementation of GGGP can be scaled up to recent clusters and supercomputers by running our experiments on multiple nodes of SuperMUC-NG, currently one of the largest supercomputing systems in Europe. 
	\item The multigrid preconditioners evolved with our method outperform all common multigrid cycles~\cite{trottenberg2000multigrid,briggs2000multigrid} with optimized relaxation factors for representative instances of the indefinite Helmholtz equation with different wavenumbers. 
	Furthermore, a subset of these methods yields a converging solver for a problem of higher difficulty and size than those considered within the optimization and for which all common multigrid cycles fail to achieve convergence.
\end{itemize} 

\section{A Formal Grammar for Generating Multigrid Preconditioners}\label{sec:grammar}
\begin{figure}
\centering
\includegraphics[width=0.85\columnwidth]{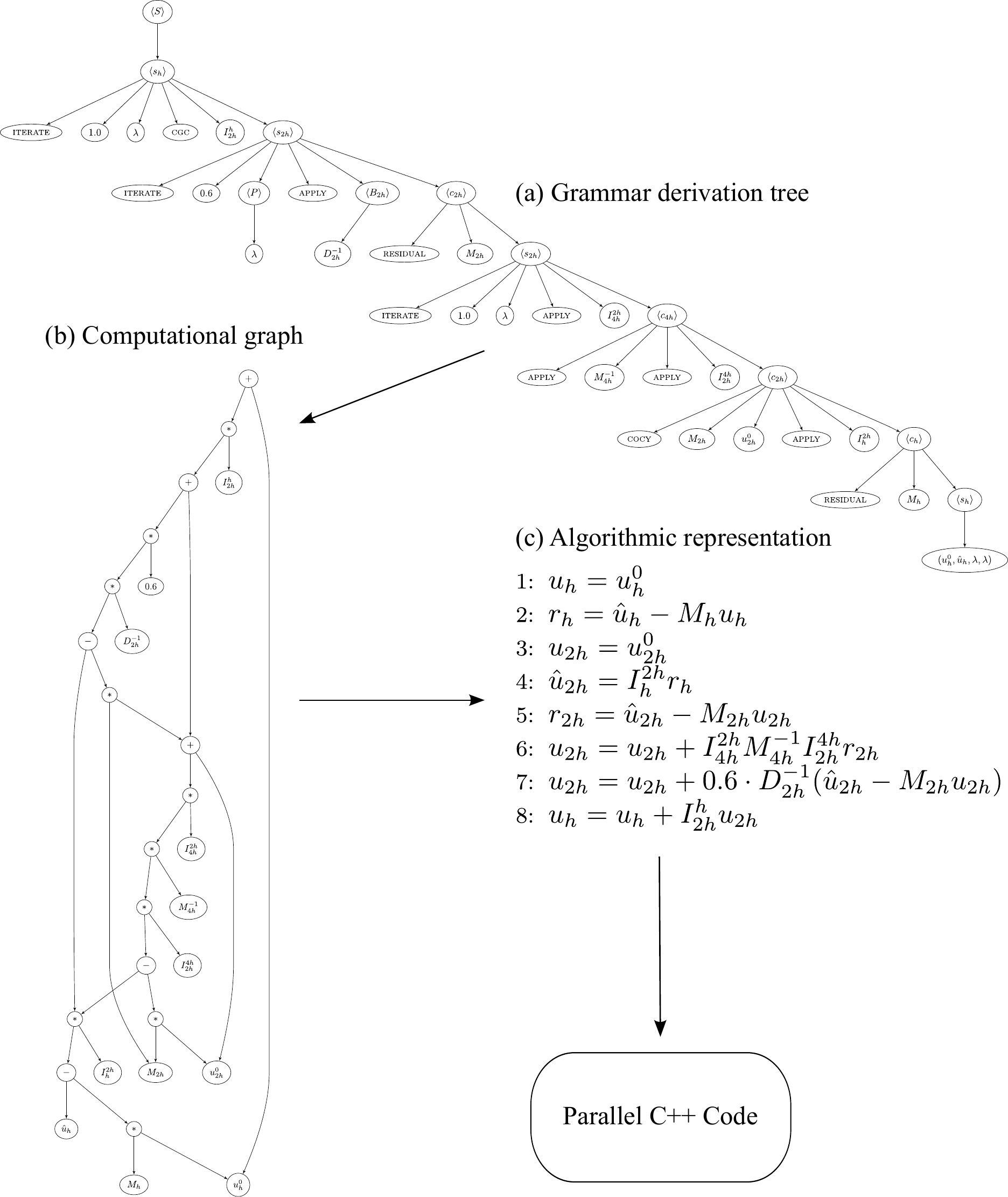}
\caption{Visualization of the process of mapping the grammar derivation tree of a three-grid V-cycle with a single step of Jacobi post-smoothing on the second finest level ($2h$) to an algorithmic representation.}
\label{fig:genotype-phenotype-mapping}
\end{figure}

We can derive a formal grammar for generating multigrid preconditioners from the one formulated in~\cite{schmitt2021evostencils,schmitt2020constructing} by replacing the system matrix $A$ with the respective preconditioning matrix $M$ and the right-hand side $f$ with $\hat{u}$. 
Table~\ref{table:grammar-syntax} contains the resulting productions for generating a multigrid preconditioner that operates on a hierarchy of three grids, where a spacing of $h$ is used on the finest grid and the only operation allowed on the coarsest grid is the application of a direct solver, denoted by the multiplication with the inverse of $M_{4h}$.
\begin{table}
	\caption{Context-free grammar for generating three-grid preconditioners.}
\label{table:grammar}
\begin{subtable}[t]{0.495\columnwidth}
	\caption{Productions}
	\begin{bnf*}
		\bnfprod{$S$} {
			\bnfpn{$s_h$}
		} \\
	    \bnfprod{$s_h$} {
	    	\bnfts{\textnormal{\textsc{iterate}}}(\bnfts{$\omega$}, \bnfpn{$P$}, \bnfts{\textnormal{\textsc{apply}}}(\bnfpn{$B_h$}, \bnfpn{$c_h$})) \bnfor
    	} \\
        \bnfmore {
        	\bnfts{\textnormal{\textsc{iterate}}}(\bnfts{$\omega$}, \bnfes, \bnfts{\textnormal{\textsc{cgc}}}(\bnfts{$I_{2h}^h$}, \bnfpn{$s_{2h}$})) \bnfor (\bnfts{$u_h^0$}, \bnfts{$\hat{u}_h$},\bnfes, \bnfes)
        } \\
		\bnfprod{$c_h$} {
			\bnfts{\textnormal{\textsc{residual}}}(\bnfts{$M_h$}, \bnfpn{$s_h$})
		} \\
		\bnfprod{$B_h$} {
			\bnfts{\textnormal{\textsc{inverse}}}(\bnfts{$M_h^{+}$}) \bnfsp \bnfts{\textnormal{with}} \bnfsp \bnfts{$M_{h} = M_{h}^{+} + M_{h}^{-}$}
		} \\
		\bnfprod{$c_{2h}$} {
			\bnfts{\textnormal{\textsc{residual}}}(\bnfts{$M_{2h}$}, \bnfpn{$s_{2h}$}) \bnfor
		} \\
		\bnfmore {
			\bnfts{\textnormal{\textsc{cocy}}}(\bnfts{$M_{2h}$},\bnfts{$u^{0}_{2h}$},\bnfts{\textnormal{\textsc{apply}}}(\bnfts{$I_h^{2h}$}, \bnfpn{$c_h$}))
		} \\
		\bnfprod{$s_{2h}$} {
			\bnfts{\textnormal{\textsc{iterate}}}(\bnfts{$\omega$}, \bnfpn{$P$}, \bnfts{\textnormal{\textsc{apply}}}(\bnfpn{$B_{2h}$}, \bnfpn{$c_{2h}$})) \bnfor
		} \\
		\bnfmore {
			\bnfts{\textnormal{\textsc{iterate}}}(\bnfts{$\omega$}, \bnfes, \bnfts{\textnormal{\textsc{apply}}}(\bnfts{$I_{4h}^{2h}$}, \bnfpn{$c_{4h}$}))
		} \\
		\bnfprod{$B_{2h}$} {
			\bnfts{\textnormal{\textsc{inverse}}}(\bnfts{$M_{2h}^{+}$}) \bnfsp \bnfts{\textnormal{with}} \bnfsp \bnfts{$M_{2h} = M_{2h}^{+} + M_{2h}^{-}$}
		} \\
		\bnfprod{${c}_{4h}$} {
			\bnfts{\textnormal{\textsc{apply}}}(\bnfts{$M^{-1}_{4h}$}, \bnfts{\textnormal{\textsc{apply}}}(\bnfts{$I_{2h}^{4h}$}, \bnfpn{$c_{2h}$}))
		} \\
		\bnfprod{$P$} {
			\bnfts{\textnormal{\textsc{partitioning}}} \bnfor \bnfes
		}
	\end{bnf*}
	\label{table:grammar-syntax}
\end{subtable}
\begin{subtable}[t]{0.495\columnwidth}
	\caption{Semantics}
	\begin{algorithmic}
		\Function{iterate}{$\omega$, $P$, $(u, \hat{u}, \delta, state)$}
		\State $\tilde{u} \gets u + \omega \cdot \delta$ with $P$
		\State return $(\tilde{u}, \hat{u}, \lambda, state)$
		\EndFunction
		\Function{apply}{$B$, $(u, \hat{u}, \delta, state)$}
		\State $\tilde{\delta} \gets B \cdot \delta$  
		\State return ($u$, $\hat{u}$, $\tilde{\delta}$, $state$)
		\EndFunction
		\Function{residual}{$M$, $(u, \hat{u}, \lambda, state)$}
		\State $\delta \gets \hat{u} - M u$
		\State return $(u, \hat{u}, \delta, state)$
		\EndFunction
		\Function{cocy}{$M_{H}$, $u_H^0$, $(u_h, \hat{u}_{h}, \delta_H, state_h)$}
		\State $u_H \gets u_H^0$
		\State $\hat{u}_{H} \gets \delta_H$
		\State $\tilde{\delta}_H \gets \hat{u}_{H} - M_H u_H$ 
		\State $state_H \gets$ $(u_h, \hat{u}_{h}, \lambda, state_h)$
		\State return $(u_H, \hat{u}_{H}, \tilde{\delta}_H, state_H)$
		\EndFunction
		\Function{cgc}{$I_H^{h}$, $(u_H, \hat{u}_{H}, \lambda, state_H)$}
		\State $(u_h, \hat{u}_{h}, \lambda, state_h) \gets state_H$
		\State $\delta_h \gets I_H^{h} \cdot u_H$
		\State return $(u_h, \hat{u}_{h}, \delta_h, state_h)$
		\EndFunction
	\end{algorithmic}
	\label{table:grammar-semantics}
\end{subtable}
\end{table}
Each rule then defines the set of expressions by which a certain variable, denoted by $\langle \cdot \rangle$, can be replaced.
Starting with the symbol $\langle S \rangle$, each expression can be substituted recursively according to the specified rules until it contains either exclusively terminal symbols or the empty string $\lambda$~\cite{linz2006introduction}.
The resulting derivation tree uniquely represents a multigrid preconditioner on the specified hierarchy of grids.
To obtain a grammar for generating multigrid preconditioners that operate on an extended hierarchy, for instance a four or five-grid method, we have to replicate the production rules formulated on the second finest level ($2h$) in Table~\ref{table:grammar-syntax} for each subsequent one.
Similar as in~\cite{schmitt2020constructing,schmitt2021evostencils} we can then formulate semantic evaluation rules, which are shown in Table~\ref{table:grammar-semantics}.
These rules guide the derivation of the corresponding sequence of multigrid operations obtained in form of a directed acyclic graph (DAG).
Figure~\ref{fig:genotype-phenotype-mapping} illustrates the resulting process of algorithm generation with the example of a three-grid V-cycle that performs a single underrelaxed Jacobi post-smoothing step on the second finest discretization level.
The resulting algorithm is formulated in Figure~\ref{fig:genotype-phenotype-mapping}c, while Figure~\ref{fig:genotype-phenotype-mapping}a shows the corresponding derivation tree based on the productions formulated in Table~\ref{table:grammar-syntax}.
Through recursive application of the rules in Table~\ref{table:grammar-semantics} the computational DAG shown in Figure~\ref{fig:genotype-phenotype-mapping}b is obtained.
Note that the recursive application of these rules in the end always returns a tuple of the form $(u_h, \hat{u}_h, \lambda, \lambda)$, whereas $u_h$ is the resulting computational graph.
Based on this intermediate representation, we can obtain an algorithmic formulation of the corresponding multigrid preconditioner by introducing variables for the approximate solution $u$, right-hand side $\hat{u}$ and residual $r$ on each discretization level, which, again, leads to Figure~\ref{fig:genotype-phenotype-mapping}c.
By making use of recent code generation techniques, it is then possible to automate the generation of optimized C++ code for a given multigrid-based solver specified in an algorithm-like fashion~\cite{kostler2020code,lengauer2020exastencils}.

Finally, note that while Table~\ref{table:grammar-syntax} fixes the number of coarsening steps until the respective problem can be solved directly, all operations are formulated relative to the discretization width $h$.
As a consequence, the computational structure of the resulting preconditioner is independent of the actual size of the grid. 
It is, therefore, possible to translate a multigrid preconditioner formulated on a hierarchy of grids with a certain depth to another one consisting of different-sized grids of the same depth.
For this purpose, we only have to replace the initial approximate solution $u_h^0$, operator $M_h$ and right-hand side $\hat{u}_h$ in Table~\ref{table:grammar-syntax} by their counterparts and reformulate each operator on the respective grid within the new hierarchy.
For instance, the derivation tree in Figure~\ref{fig:genotype-phenotype-mapping}a can be translated to the computational DAG of a structurally similar multigrid preconditioner for a different problem discretized on a hierarchy of three grids. 
Consequently, every multigrid method produced by a grammar formulated on a particular hierarchy of grids is generalizable over the set of all structurally equivalent grammars that employ the same number of coarsening steps.
We can thus apply a multigrid preconditioner obtained on a specific instance of the Helmholtz equation to similar problems. 
In the following, we utilize this principle to evolve efficient multigrid-based preconditioners that can be generalized to a whole class of Helmholtz problems.
\section{Evolutionary Search Method}
\label{sec:evolutionary-search}
After establishing a generalizable representation for multigrid preconditioners of arbitrary structure, the next task is to formulate a search method that can identify those leading to an efficient solver for different Helmholtz problem instances.
While the branching factor of the productions formulated in Table~\ref{table:grammar-syntax} may seem small at first, it has already been shown in~\cite{schmitt2021evostencils} that, in practice, the resulting search space exceeds any size for which a simple exhaustive search is applicable.
In~\cite{schmitt2021evostencils} we have provided $3 \cdot 10^{14}$ as a lower bound for the size of the search space of a three-grid method with a limited number of choices for smoothing on each level.
Since the construction of a multigrid preconditioner comprises a similar number of choices, this lower bound also applies in the present case, rendering a mere brute-force search unfeasible.
Search heuristics can often find an acceptable approximation for the global optimum when the search space is too large to evaluate all possible solutions.
In principle, the quality of a preconditioner can be assessed by considering two objectives. 
First of all, preconditioning aims to minimize the condition number of the matrix $A M^{-1}$.
Since, in practice, the inverse $M^{-1}$ is not computed explicitly, the effectiveness of a preconditioner depends on its approximation accuracy, which directly affects the \emph{number of iterations} required by an iterative method to achieve a certain error reduction.
On the other hand, a method that achieves the same quality of approximation but can be executed faster on modern computer architectures achieves a lower \emph{execution time per iteration}.
The task of finding an optimal multigrid preconditioner can thus be considered as a multi-objective search problem.
\paragraph{Fitness Evaluation and Generalization}
As shown in~\cite{erlangga2008advances} the number of iterations required for solving Equation~\eqref{eq:helmholtz} using a preconditioned solver grows with the wavenumber $k$.
This work aims to obtain multigrid methods that can be generalized over different problem instances. 
Therefore, before evaluating a given preconditioner on problems with a large wavenumber, we first consider an instance of the same problem with a smaller wavenumber and hence lower difficulty.
When we start with a random initialization, the probability of generating multigrid methods that do not represent effective preconditioners is high, even for a problem instance that is, in principle, comparably easy to solve.
As the search progresses, the average quality of the obtained preconditioners is expected to improve. 
Hence, the probability that they can also be successfully applied to instances with higher difficulty increases.
On the other hand, most multigrid methods that are efficient in preconditioning a problem instance with a large wavenumber can be expected to function also on a problem with a smaller wavenumber. 
We, therefore, propose a stepwise adaption of the difficulty of the evaluated problem.
To perform the actual search, we employ a multi-objective GGGP-based algorithm that operates on a population of derivation trees~\cite{mckay2010grammar,whigham1995grammatically}.
Each tree represents a certain point in the search space considered at the current step of the method.
In every new step, the search progresses by creating a new population of trees based on the current one.
The resulting procedure is summarized in Algorithm~\ref{alg:evolutionary-search}, in which the difficulty of the problem considered for evaluation is adjusted in every $m$th iteration of the search.
\begin{algorithm}
	\caption{Evolutionary Search}
	\label{alg:evolutionary-search}
	\begin{algorithmic} 
		\State \textbf{Construct} the grammar $G_0$ for the initial problem
		\State \textbf{Initialize} the population $P_0$ based on $G_0$
		\State \textbf{Evaluate} $P_0$ on the initial problem
		\For{$i := 0, \dots, n$}
		\If{$i > 0$ and $i \mod m = 0$}
		\State $j := i / m$ 
		\State Increase the problem difficulty
		\State \textbf{Construct} the corresponding grammar $G_j$
		\State  \textbf{Adapt} the current population $P_i$ to $G_j$
		\State \textbf{Evaluate} $P_i$ on the new problem
		\EndIf
		\State \textbf{Generate} new solutions $C_i$ based on $P_i$ and $G_j$
		\State \textbf{Evaluate} $C_i$ on the current problem
		\State \textbf{Select} $P_{i+1}$ from $C_i \cup P_i$
		\EndFor
	\end{algorithmic}
\end{algorithm}

The question that remains to be answered is how a new population $P_{i+1}$ should be generated in each step of Algorithm~\ref{alg:evolutionary-search} based on the current one.
In principle, the current population $P_i$ represents the subspace of possible multigrid preconditioners considered within step $i$ of the search.
Accordingly, the generation of $P_{i+1}$ represents moving the search to a new subspace, which is expected to contain solutions that, according to both objectives, correspond to multigrid methods representing more efficient preconditioners for the given problem than those located in the current subspace.
Consequently, we need to evaluate the quality of a subspace represented by the current population in terms of its potential to obtain efficient preconditioners from it.
In principle, the number of iterations required to achieve a particular error reduction with a preconditioned iterative method could be predicted with local Fourier analysis~\cite{cools2013analysis}. 
However, there has been only a limited amount of research on the accuracy of this method for evaluating nonstandard multigrid methods.
In particular, the experiments performed in~\cite{schmitt2020constructing} indicate that the predictions obtained with this method are not always consistent with experimentally determined behavior.
Alternatively, we can obtain all relevant performance characteristics of a solver through its direct application to a representative test problem.
For this purpose, it is necessary to automatically generate an implementation based on the algorithmic representation of those multigrid preconditioners obtained through semantical evaluation of the derivation trees produced by the respective grammar.
ExaStencils~\cite{lengauer2020exastencils} is a framework that has been specifically designed for the automatic generation of scalable multigrid implementations based on a tailored domain-specific language (DSL) called ExaSlang~\cite{schmitt2014exaslang}.
It enables the specification of solvers in a discretization level-independent manner. 
At the same time, the actual size of the problem can be controlled utilizing simple configuration files, which grants us the possibility to automatically generate implementations for a specific solver that are executable on a wide range of different computer architectures.
These can then be evaluated on the target platform for the two objectives, i.e., number of iterations and execution time per iteration.
\paragraph{Implementation of Grammar-Guided Genetic Programming}
In each iteration of the search, we create a new population based on the existing one using GGGP, where each individual represents a derivation tree of the form of Figure~\ref{fig:genotype-phenotype-mapping}a.
To apply this method to the grammar formulated in Table~\ref{table:grammar-syntax}, we first need to consider its unique structure.
Note that except for the variables $\langle s_h \rangle$, $\langle B_h \rangle$, $\langle B_{2h} \rangle$ and $\langle P \rangle$ none of the expressions generated by any of the productions of a variable exclusively consists of terminals.
Consequently, the grammar does not permit the construction of a derivation tree with branches of equal length.
We, therefore, employ the \emph{grow} strategy as described in~\cite{poli2008field,koza1992genetic} to initialize the population.
A derivation tree is extended by randomly choosing a production from the combined set of terminal and non-terminal productions until the longest path within the tree exceeds a certain depth.
Since the grammar comprises an even branching factor of two for all non-terminal productions, there is no need to adapt the probability of selecting a particular production, but choosing uniformly from the combined set already results in a sufficient diversity in the population.

To create new individuals based on an existing population, we employ mutation and recombination.
For this purpose, we first select several individuals using a binary tournament selection based on the dominance relation and crowding distance between individuals, as described in~\cite{deb2002fast}.
Mutation is performed by randomly selecting a variable node within the given derivation tree.
The subtree for which this variable represents the root node is then replaced by a new randomly generated tree, which is created using the \emph{grow} initialization operator. 
However, we also permit the insertion of the replaced subtree as a branch within the new one.
Note that this insertion, which is only allowed once, is only possible if the variable that represents the root node of the original subtree occurs within the new one.
Consequently, if this condition is never fulfilled, the original subtree is replaced without insertion.
Therefore, our mutation operator can either perform subtree replacement or insertion in case the newly generated subtree can connect the original one to its root node. 
While mutation is performed on a single individual, within recombination we create two new individuals by combining the derivation trees of two individuals selected as \emph{parents}.
For this purpose, we employ standard subtree crossover as described in~\cite{poli2008field}, whereby we choose the crossover point uniformly among all possible nodes within both trees.
Finally, after the creation and evaluation of a certain number of novel solutions, we employ the sorting procedure described in~\cite{fortin2013generalizing} to identify the non-dominated solutions in the combined set of the newly created and existing ones. 
These individuals then form the new population $P_{i+1}$ in the next step of Algorithm~\ref{alg:evolutionary-search}.

\section{Experimental Evaluation}
To evaluate the effectiveness of our evolutionary search method in finding multigrid methods that act as efficient preconditioners, we consider the two-dimensional Helmholtz equation on a unit square with Dirichlet boundary conditions at the top and bottom, and Robin radiation conditions at the left and right, as defined by
\begin{equation*}
	\label{eq:helmholtz-test-problem}
	\begin{split}
		(-\nabla ^{2} - k^{2}) u & = f \quad \text{in} \; \left( 0, 1 \right)^2 \\
		u & = 0 \quad \text{on} \; \left( 0, 1 \right) \times \{0\}, \left( 0, 1 \right) \times \{1\} \\
		\partial_{\mathbf{n}} u - iku & = 0 \quad \text{on} \; \{0\} \times \left( 0, 1 \right), \{1\} \times \left( 0, 1 \right) \\
		f(x, y) & = \delta(x - 0.5, y - 0.5),
	\end{split}
\end{equation*}
where $\delta(\mathbf{x})$ represents the Dirac delta function.
We discretize this equation on a uniform Cartesian grid using the classical five-point stencil
\begin{equation*}
	\frac{1}{h^2} \begin{bmatrix}
		& -1 & \\
		-1 & 4 - (k h)^2 & -1 \\
		& -1 &  
	\end{bmatrix},
\end{equation*}
while the Dirac delta function is approximated with a second-order Zenger correction~\cite{koestler2004extrapolation}.
The step size $h$ of the grid is chosen to fulfill the second-order accuracy requirement $h k = 0.625$ as suggested in~\cite{erlangga2006multigrid}. 
Since the analytic solution of this equation is not known in advance, we consider an approximate solution to be sufficient if the initial residual has been reduced by a factor of $10^{-7}$ for $k  \leq 160$ and $10^{-6}$ for all larger wavenumbers.
The resulting complex-valued system of linear equations is indefinite, and the required number of iterations for solving it with a non-preconditioned Krylov subspace method increases drastically with the wavenumber $k$~\cite{erlangga2006multigrid,erlangga2008advances}.  
As a solver, we, therefore, employ a biconjugate gradient stabilized method (BiCGSTAB)~\cite{saad2003iterative}, right-preconditioned with a shifted Laplacian
\begin{equation*}
	M = -\nabla ^{2} - (k^{2} + 0.5 i k^{2}),
\end{equation*}
which is among the suggested solvers in~\cite{erlangga2008advances}.
The resulting numerical solution method is summarized in Algorithm~\ref{alg:preconditioned-bicgstab}.
Note that greek letters denote scalar values while small letters represent vectors and capital letters matrices.
\begin{algorithm}[t]
	\caption{Right-Preconditioned BiCGSTAB}
	\label{alg:preconditioned-bicgstab}
	\begin{algorithmic}[1] 
		\State $\hat{r}^0 = r^0 = f - A u^0$
		\State $\alpha_0 = \beta_0 = \rho_0 = \omega_0 = 1$
		\State $p^0 = q^0 = 0$
		\For{$i := 1, \dots, n$}
		\State $\rho_i = \hat{r}^{i-1} \cdot r^{i-1}$
		\State $\beta_i = \frac{\rho_i }{\rho_{i-1} }\frac{\alpha_{i-1}}{ \omega_{i-1}}$
		\State $p^i = r^{i-1} + \beta_i (p^{i-1} - \omega_{i-1} q^{i-1})$
		\State Solve $M u^i = p^i$
		\State $q^i = A u^i$
		\State $\alpha_i = \rho_i / (\hat{r}^{i-1} \cdot q^i)$
		\State $h^i = u^{i-1} + \alpha_i u^i$	
		\State $s^i = r^{i-1} - \alpha_i r^i$
		\State Solve $M u^i = s^i$
		\State $t^i = A u^i$
		\State $\omega_i = (t^i \cdot s^i) / (t^i \cdot t^i)$
		\State $u^i = h^i + \omega_i u^i$
		\State $r^i = s^i - \omega_i t^i$
		\If{$\norm{r^i}/\norm{r^0} < \varepsilon$}
		\Return $u^i$
		\EndIf
		\EndFor
	\end{algorithmic}
\end{algorithm}
In each step of this iterative scheme, it is necessary to compute an approximate solution for two systems of linear equations of the form $M u = \hat{u}$, in line 8 and 13, 
each of which is achieved through the application of a single multigrid iteration.
\subsection{Optimization Settings}
To evaluate the behavior of our GGGP-based search method, we perform a total number of ten randomized optimization runs.
While we are aware that this number is insufficient for a reasonable statistical evaluation of an evolutionary algorithm's behavior, it still enables us to demonstrate that our method is capable of repeatedly evolving generalizable preconditioners for the given test problem.
Even though a more accurate assessment of our method's behavior would be desirable, the high computational and temporal costs of each run, which each take between 24 and 48 hours, put a strict limit on the number of experiments.
Within all optimization runs, we choose a step size of $h = 1/2^{l}$ on each level $l$, whereby we employ a range of $l \in \left[l_{max} - 4,l_{max}\right]$. 
Accordingly, our goal is to construct an optimal five-grid preconditioner for the given problem.
For this purpose, we consider the following components:
\begin{description}
	\item[\textbf{Smoothers}:] Pointwise and block Jacobi with rectangular blocks up to a maximum number of six terms, red-black Gauss-Seidel
	\item[\textbf{Restriction}:] Full-weighting restriction
	\item[\textbf{Prolongation}:] Bilinear interpolation
	\item[\textbf{Relaxation factors}:] $\left( 0.1 + 0.05i \right)_{i = 0}^{36} = \left(0.1, 0.15, 0.2, \dots, 1.9 \right)$
	\item[\textbf{Coarse-grid solver}:] BiCGSTAB for $l = l_{max} - 4$
\end{description}
To generate block Jacobi smoothers, we define a splitting $M_h = L_h + D_h + U_h$ where $D_h$ is a block diagonal matrix, such that we have to solve a local system whose size corresponds to the size of a block at every grid point.
For a more detailed treatment of block relaxation methods, the reader is referred to~\cite{trottenberg2000multigrid,saad2003iterative}.
The relaxation factor $\omega$ for each smoothing and coarse-grid correction step is chosen from the above sequence.
By employing the same code generation-based optimizations for each component of a solver, we ensure that the resulting measurements are comparable for all multigrid variants considered in this work.
All experiments are performed on the SuperMUC-NG cluster
, where each node represents an Intel Skylake Xeon Platinum 8174 processor that consists of eight islands, each with six physical cores.
While within the optimization, we evaluate each individual's fitness on a single island, the final evaluation is performed on a full node of the system with 48 cores.
We employ GCC 7.5 as a compiler, using the -O3 optimization level and an OpenMP-based parallelization in both cases. 
To assess each preconditioner's generalizability, we consider the three different wavenumbers 160, 320, and 640, together with a discretization width of $h = 0.625 / k$.
\subsection{Reference Methods}
To establish a baseline, we consider several well-known and commonly used multigrid cycles~\cite{trottenberg2000multigrid} that are all based on the application of a certain smoother for a fixed number of times.
The resulting solver is translated to ExaStencils' DSL, based on which a multi-threaded C++ implementation is generated and executed on a full SuperMUC-NG node using 48 OpenMP threads. 
To optimize each multigrid cycle's effectiveness as a preconditioner, we experimentally obtain the optimum relaxation factor for $k = 320$ from the mentioned interval.
While in~\cite{erlangga2008advances,cocquet2017shift} a damped Jacobi is employed as a smoother, in the given case, it does not result in a convergent solver for $k > 80$. 
We have verified this assumption by considering every possible relaxation factor value from the given interval.
In contrast, red-black Gauss-Seidel represents an effective smoother for the considered range of $k$.
The second column of Table~\ref{table:reference-methods} contains the optimum red-black Gauss-Seidel relaxation factor ($\omega$) for each cycle. 
For instance, V(2, 1) represents a V-cycle that performs two pre- and one post-smoothing step on each discretization level.
Using the same relaxation factor, we employ each cycle as a preconditioner for the three wavenumbers considered.
For consistent measurements, each solver is executed ten times to compute the average solving time of all runs, which reduces the deviations to a negligible level.
The results are shown in the remaining columns of Table~\ref{table:reference-methods}.
Here omitted values imply that the corresponding preconditioned BiCGSTAB method did not achieve the required error reduction within 20,000 iterations.
\begin{table}
	\caption{Comparison of multigrid preconditioners - Number of iterations and average time required for solving a problem with the particular wavenumber.}
	\begin{subtable}[t]{0.445\columnwidth}
		\caption{Reference methods with optimum relaxation factors $\omega$ for $k = 320$}
		\label{table:reference-methods}
		\centering
		\begin{tabular}{l l l l l l}
			\toprule
			& $\omega$ & \multicolumn{2}{c}{Iterations} & \multicolumn{2}{c}{Solving Time (s)} \\
			\cmidrule(r){3-4} \cmidrule(r){5-6}
			$k$ & & $160$ & $320$ & $160$ & $320$ \\
			\midrule
			V(0, 1) & $1.25$ & $2078$ & $6297$ & $6.38$ & $35.11$ \\
			\midrule
			V(1, 1) & $0.6$ & $1880$ & $6297$ & $7.66$ & $44.27$ \\
			\midrule
			V(2, 1) & $0.6$ & $-$ & $5532$ & $-$ & $47.0$ \\
			\midrule
			V(2, 2) & $0.5$ & $1627$ & $5115$ & $9.93$ & $50.54$  \\
			\midrule
			V(3, 3) & $0.4$ & $1753$ & $5168$ & $13.97$ & $76.00$ \\
			\midrule
			F(0, 1) & $1.15$ & $1467$ & $4028$ & $8.15$ & $42.87$  \\
			\midrule
			F(1, 1) & $0.75$ & $1546$ & $3988$ & $11.21$ & $54.51$ \\
			\midrule
			F(2, 1) & $0.55$ & $1146$ & $3934$ & $10.87$ & $67.62$ \\
			\midrule
			F(2, 2) & $0.65$ & $1060$ & $3213$ & $13.92$ & $65.06$ \\
			\midrule
			F(3, 3) & $0.45$ & $1085$ & $3464$ & $18.88$ & $92.97$ \\
			\midrule
			W(0, 1) & $0.75$ & $1265$ & $4215$ & $8.67$ & $72.08$ \\
			\midrule
			W(1, 1) & $0.8$ & $1208$ & $3570$ & $13.08$ & $76.22$ \\
			\midrule
			W(2, 1) & $0.6$ & $1313$ & $3074$ & $17.71$ & $79.67$ \\
			\midrule
			W(2, 2) & $0.5$ & $1069$ & $3376$ & $17.14$ & $101.6$ \\
			\midrule
			W(3, 3) & $0.45$ & $942$ & $2976$ & $19.65$ & $117.8$ \\
			\bottomrule
		\end{tabular}
	\end{subtable}
\begin{subtable}[t]{0.545\columnwidth}
		\caption{Best preconditioners according to the product of both objectives}
	\label{table:evolved-solvers}
	\centering
	\begin{tabular}{l l l l l l l }
		\toprule
		& \multicolumn{3}{c}{Iterations} & \multicolumn{3}{c}{Solving Time (s)} \\
		\cmidrule(l){2-4} \cmidrule(l){5-7}
		$k$ & $160$ & $320$ & $640$ & $160$ & $320$ & $640$ \\
		\midrule
		EP-1 & $1178$ & $3399$ & $-$ & $6.29$ & $28.07$ & $-$ \\
		\midrule
		EP-2 & $795$ & $2160$ & $8449$ & $7.86$ & $29.89$ & $241.7$\\
		\midrule
		EP-3 & $933$ & $2827$ & $11143$ & $6.08$ & $27.58$ & $257.8$ \\
		\midrule
		EP-4 & $637$ & $2509$ & $7901$ & $7.17$& $41.04$ & $268.2$ \\
		\midrule
		EP-5 & $539$ & $1838$ & $7765$ & $5.01$ & $28.39$ & $227.7$ \\
		\midrule
		EP-6 & $941$ & $2103$ & $-$ & $9.58$ & $30.76$ & $-$ \\
		\midrule
		EP-7 & $955$ & $2701$ & $-$  & $6.45$& $27.84$ & $-$ \\
		\midrule
		EP-8 & $945$ & $2870$ & $10839$ & $7.24$ & $33.02$ & $276.9$ \\
		\midrule
		EP-9 & $3436$ & $3872$ & $-$  & $15.15$ & $27.51$ & $-$ \\
		\midrule
		EP-10 & $586$ & $1881$ & $8855$ & $6.70$ & $31.39$ & $246.1$ \\
		\bottomrule
	\end{tabular}
\end{subtable}
\end{table}
Additionally, we have evaluated each resulting cycle on a problem with wavenumber $k = 640$, which did not yield a convergent solver in any of the cases.
For the other two wavenumbers considered, the W(3, 3)-cycle represents the most effective preconditioner and therefore leads to the lowest number of iterations, while the V(0, 1)-cycle yields the overall fastest solving time on the given platform.
\subsection{Implementation Details}
We have implemented the evolutionary search procedure summarized in Algorithm~\ref{alg:evolutionary-search} by extending the approach described in~\cite{schmitt2021evostencils}.
For this purpose, we generate a new grammar for constructing complex-valued multigrid preconditioners adapted to the respective problem instance in every $m$th iteration.
For evaluating each method's fitness, it is first translated to an ExaSlang~\cite{schmitt2014exaslang} representation, which is then automatically integrated into an existing Krylov Subspace method as a preconditioner.
Finally, based on the DSL representation of each solver, a multi-threaded C++ implementation is generated, as described in Section~\ref{sec:evolutionary-search}.
To cope with the cost of running a code generation pipeline for each evaluation, together with the growing execution time required for solving the increasingly difficult problem instances, we employ a distributed parallelization using the Message Passing Interface (MPI) library.
The resulting implementation is freely available as part of the open-source library EvoStencils\footnote{EvoStencils: \url{https://github.com/jonas-schmitt/evostencils}}.

For each experiment, we perform an evolutionary search with a total population size of 128 on eight nodes of SuperMUC-NG using 64 MPI processes, whereby each process is executed on a separate island of the system.
As an initial problem, we choose $k = 80$ with a maximum level $l_{max} = 7$, discretized with a step size $h = 1/2^{7}$.
A problem instance with greater difficulty is then constructed by doubling the wavenumber.
Note that due to the requirement $hk = 0.625$, this results in a step size half as large as the original one and, in total, a four times larger grid.
The relaxation factor for each smoothing and coarse-grid correction step is chosen from the interval specified above.
Each preconditioner is evaluated on the respective island using 12 threads.
A solver is considered convergent for all wavenumbers if it can reduce the initial residual by $10^{-7}$ in less than 10,000 iterations. 
The initial population is obtained through the non-dominated sorting of a randomly generated set of 1024 individuals.
The search is then performed for 150 iterations, whereby the difficulty is adjusted every 50 iterations.
New derivation trees are created through recombination with a probability of $2/3$ or by mutation.
In the latter, a terminal symbol is chosen with a probability of $1/3$.
To evaluate the consistency of the obtained results, we perform ten experiments with a random initialization.
At the end of each experiment, we identify the ten best preconditioners according to the product of both objectives and evaluate them under the same conditions as the reference methods, i.e., by executing each solver ten times on a full SuperMUC-NG node using 48 OpenMP threads.

\section{Results and Discussion}
Table~\ref{table:evolved-solvers} contains the results for each preconditioner from the set of non-dominated solutions of the respective experiment that achieves the fastest solving time for a wavenumber of $k = 640$.
\begin{figure}
\centering
	\includegraphics[scale=0.75]{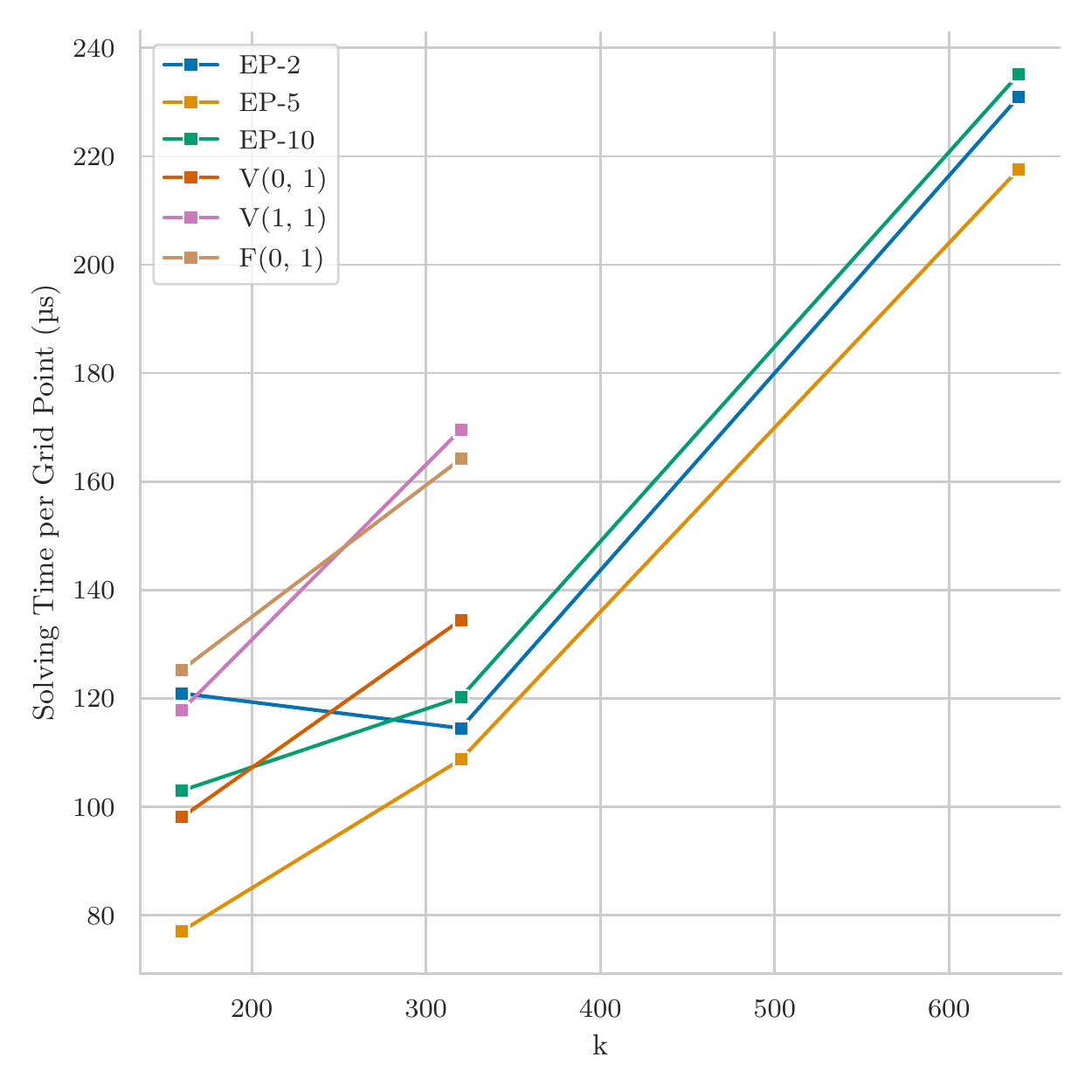}
	\caption{Solving time comparison of the best preconditioners according to the product of both objectives for different wavenumbers ($k$).}
	\label{fig:solving-time}
\end{figure}
Note that, similar to all reference methods, in four of the ten cases, EP-1, EP-6, EP-7, and EP-9, none of the evaluated solvers could achieve convergence for the largest wavenumber $k = 640$.
In these cases, we have selected the preconditioner, which leads to the fastest solving time for $k = 320$.
Additionally, Figure~\ref{fig:solving-time} shows a direct comparison of the best-performing preconditioners from both groups for different wavenumbers.
For better comparability of the results achieved on different grid sizes, we measure the solving time per grid point instead of the total time required for solving each problem.
Note that all solvers have been evaluated using the same number of measurements on the target platform.
All three evolved preconditioners included in this plot represent more efficient methods for $k \geq 320$ than the best of the reference methods, the V(0, 1)-cycle, while remaining competitive for lower wavenumbers.
The most efficient preconditioner (EP-5) leads to a consistent improvement of about 20 \% compared to the V(0, 1)-cycle for all wavenumbers.
Furthermore, while none of the evolved preconditioners has been evaluated on wavenumbers greater than 320 within the search, in six of the ten experiments they lead to a converging solution method for the case of a wavenumber of $640$, for which all standard methods fail.
In the remaining four experiments, the search still finds competitive preconditioners that generalize well for $k \leq 320$, whereby only in one case (EP-9) preconditioning leads to an inefficient solver for a wavenumber of 160.
Therefore, we have demonstrated that our evolutionary search method could find generalizable and efficient methods in the majority of the experiments performed.
In addition, our GGGP-based approach was able to evolve methods that surpass the capabilities of standard multigrid cycles in preconditioning Helmholtz problems, which could be demonstrated by solving a problem instance with $k = 640$.

\begin{figure}
	\centering
	\includegraphics[scale=0.75]{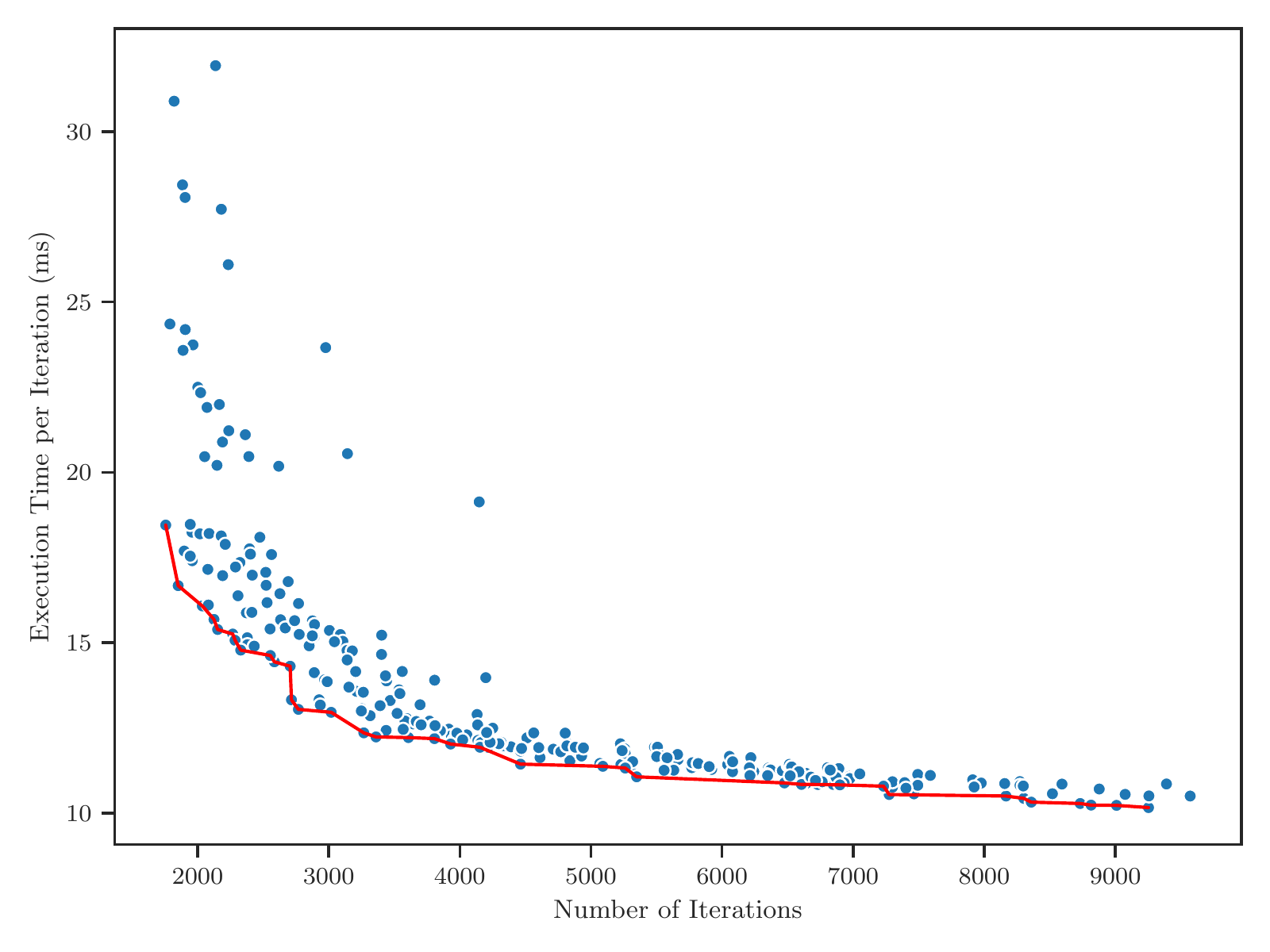}
	\caption{Distribution of non-dominated solutions at the end of all ten experiments for $k = 320$. The red line denotes the combined front.}
	\label{fig:pareto-front}
\end{figure}

To further investigate our evolutionary algorithm's behavior, Figure~\ref{fig:pareto-front} shows the distribution of the non-dominated solutions at the end of all ten experiments, whereby the red line represents the non-dominated front of the combined set of solutions.
From this distribution, we can conclude that the method can consistently find preconditioners of equal quality for iteration numbers of 4000 or more.
For lower iteration numbers, the spread between the solutions increases and, in individual experiments, leads to suboptimal preconditioners in the left part of the solution space.
Reducing the number of iterations requires more smoothing and coarse-grid corrections steps with an effective combination of relaxation factors, which results in a growth of the size of the corresponding grammar derivation tree.
This effect impedes the search algorithm's ability to find the right combination of productions that leads to the same Pareto-optimal preconditioner in every experiment. 

While the efficiency and generalizability of the evolved preconditioners could be demonstrated, we have not yet analyzed how these methods function algorithmically.
For this purpose, we consider the two evolved preconditioners, EP-5 and EP-2, that perform best for wavenumbers $k \geq 320$.
\begin{figure}
	\begin{subfigure}[t]{0.85\columnwidth}
		\centering
		\resizebox{\columnwidth}{!}{%
			\begin{tikzpicture}
				\node   (h) at (-0.75, 4){$h$};
				\node   (2h) at (-0.75, 3){$2h$};
				\node   (4h) at (-0.75, 2){$4h$};
				\node   (8h) at (-0.75, 1){$8h$};
				\node   (16h) at (-0.75, 0){$16h$};
				\node	(a) at (0,4) [draw, circle,inner sep=0pt,minimum size=5mm] {\phantom{\scriptsize 1.00}};
				\node	(b) at (0.5,3) [draw, circle,fill=lightred,inner sep=0pt,minimum size=5mm] {\scriptsize 1.30};
				\node	(c) at (1,2) [draw,circle,fill=lightblue, inner sep=0pt,minimum size=5mm] {\scriptsize 1.45};
				\node	(d) at (1.5,1) [draw,circle,fill=lightblue, inner sep=0pt,minimum size=5mm] {\scriptsize 1.35};
				\node	(e) at (2,0) [draw, circle,fill=black, inner sep=0pt,minimum size=5mm] {\phantom{\scriptsize 1.00}};
				\node	(f) at (2.5,1) [draw, circle, inner sep=0pt,minimum size=5mm] {\phantom{\scriptsize 1.00}};
				\node	(g) at (3,2) [draw, circle, inner sep=0pt,minimum size=5mm] {\phantom{\scriptsize 1.00}};
				\node	(h) at (3.5,3) [draw, circle, inner sep=0pt,minimum size=5mm] {\phantom{\scriptsize 1.00}};
				\node	(i) at (4,4) [draw, circle,fill=lightred, inner sep=0pt,minimum size=5mm] {\scriptsize 1.20};
				\node	(j) at (4.5,3) [draw, circle, inner sep=0pt,minimum size=5mm] {\phantom{\scriptsize 1.00}};
				\node	(k) at (5,2) [draw, circle,fill=lightred, inner sep=0pt,minimum size=5mm] {\scriptsize 1.30};
				\node	(l) at (5.5,3) [draw, circle, inner sep=0pt,minimum size=5mm] {\phantom{\scriptsize 1.00}};
				\node	(m) at (6,2) [draw, circle,fill=lightblue, inner sep=0pt,minimum size=5mm] {\scriptsize 1.80};
				\node	(n) at (7,2) [draw, circle,fill=lightblue, inner sep=0pt,minimum size=5mm] {\scriptsize 0.55};
				\node	(o) at (7.5,3) [draw, circle, inner sep=0pt,minimum size=5mm] {\phantom{\scriptsize 1.00}};
				\node	(p) at (8,2) [draw, circle, fill=lightred, inner sep=0pt,minimum size=5mm] {\scriptsize 1.80};
				\node	(q) at (8.5,1) [draw, circle, fill=lightred, inner sep=0pt,minimum size=5mm] {\scriptsize 1.45};
				\node	(r) at (9,2) [draw, circle, inner sep=0pt,minimum size=5mm] {\phantom{\scriptsize 1.00}};
				\node	(s) at (9.5,3) [draw, circle, inner sep=0pt,minimum size=5mm] {\phantom{\scriptsize 1.00}};
				\node	(t) at (10,2) [draw, circle, fill=lightblue, inner sep=0pt,minimum size=5mm] {\scriptsize 0.95};
				\node	(u) at (11,2) [draw, circle, fill=lightred, inner sep=0pt,minimum size=5mm] {\scriptsize 0.70};
				\node	(v) at (11.5,3) [draw, circle, fill=lightred, inner sep=0pt,minimum size=5mm] {\scriptsize 1.00};
				\node	(w) at (12,4) [draw, circle, fill=lightred, inner sep=0pt,minimum size=5mm] {\scriptsize 0.90};
				\draw 
				(a) edge[->] (b) 
				(b) edge[->] (c)
				(c) edge[->] (d)
				(d) edge[->] (e)   
				(e) edge[->] node[near end,left] {\scriptsize 1.00}  (f)
				(f) edge[->] node[near end,left] {\scriptsize 1.00} (g)
				(g) edge[->] node[near end,left] {\scriptsize 0.80} (h)
				(h) edge[->] node[near end,left] {\scriptsize 1.60} (i)
				(i) edge[->] (j) 
				(j) edge[->] (k)
				(k) edge[->] node[near end,left] {\scriptsize 0.85} (l)
				(l) edge[->] (m)   
				(m) edge[->] (n)
				(n) edge[->] node[near end,left] {\scriptsize 0.85} (o)
				(o) edge[->] (p)
				(p) edge[->] (q)
				(q) edge[->] node[near end,left] {\scriptsize 1.05} (r)
				(r) edge[->] node[near end,left] {\scriptsize 1.10} (s)
				(s) edge[->] (t)
				(t) edge[->] (u)
				(u) edge[->] node[near end,left] {\scriptsize 1.05} (v)
				(v) edge[->] node[near end,left] {\scriptsize 1.80} (w)
				;
			\end{tikzpicture}
		}
		\caption{EP-5}
		\label{fig:ep-5}
	\end{subfigure}
	\begin{subfigure}[t]{0.85\columnwidth}
		\centering
		\resizebox{\columnwidth}{!}{%
			\begin{tikzpicture}
				\node   (h) at (-0.75, 4){$h$};
				\node   (2h) at (-0.75, 3){$2h$};
				\node   (4h) at (-0.75, 2){$4h$};
				\node   (8h) at (-0.75, 1){$8h$};
				\node   (16h) at (-0.75, 0){$16h$};
				\node	(a) at (0,4) [draw, circle,inner sep=0pt,minimum size=5mm] {\phantom{\scriptsize 1.00}};
				\node	(b) at (0.5,3) [draw, circle,inner sep=0pt,minimum size=5mm] {\phantom{\scriptsize 1.00}};
				\node	(c) at (1,2) [draw, circle,fill=lightred,inner sep=0pt,minimum size=5mm] {\scriptsize 1.00};
				\node	(d) at (1.5,1) [draw, circle,inner sep=0pt,minimum size=5mm] {\phantom{\scriptsize 1.00}};
				\node	(e) at (2,0) [draw, circle,fill=black, inner sep=0pt,minimum size=5mm] {\phantom{\scriptsize 1.00}};
				\node	(f) at (2.5,1) [draw, circle,inner sep=0pt,minimum size=5mm] {\phantom{\scriptsize 1.00}};
				\node	(g) at (3,2) [draw, circle,inner sep=0pt,minimum size=5mm] {\phantom{\scriptsize 1.00}};
				\node	(h) at (3.5,3) [draw, circle,fill=lightred,inner sep=0pt,minimum size=5mm] {\scriptsize 0.90};
				\node	(i) at (4.5,3) [draw, circle,fill=lightred,inner sep=0pt,minimum size=5mm] {\scriptsize 1.20};
				\node	(j) at (5,2) [draw, circle,fill=lightred, inner sep=0pt,minimum size=5mm] {\scriptsize 1.40};
				\node	(k) at (5.5,1) [draw, circle,fill=lightblue,inner sep=0pt,minimum size=5mm] {\scriptsize 1.00};
				\node	(l) at (6,2) [draw, circle,inner sep=0pt,minimum size=5mm] {\phantom{\scriptsize 1.00}};
				\node	(m) at (6.5,3) [draw, circle, fill=lightred, inner sep=0pt,minimum size=5mm] {\scriptsize 1.20};
				\node	(n) at (7,2) [draw, circle,fill=lightred, inner sep=0pt,minimum size=5mm] {\scriptsize 1.40};
				\node	(o) at (7.5,1) [draw, circle, fill=lightred, inner sep=0pt,minimum size=5mm] {\scriptsize 0.65};
				\node	(p) at (8,2) [draw, circle, inner sep=0pt,minimum size=5mm] {\phantom{\scriptsize 1.00}};
				\node	(q) at (8.5,3) [draw, circle, fill=lightred, inner sep=0pt,minimum size=5mm] {\scriptsize 1.20};
				\node	(r) at (9,2) [draw, circle, fill=lightred, inner sep=0pt,minimum size=5mm] {\scriptsize 1.25};
				\node	(s) at (9.5,1) [draw, circle, fill=lightred, inner sep=0pt,minimum size=5mm] {\scriptsize 1.45};
				\node	(t) at (10,2) [draw, circle, inner sep=0pt,minimum size=5mm] {\phantom{\scriptsize 1.00}};
				\node	(u) at (10.5,3) [draw, circle, fill=lightred, inner sep=0pt,minimum size=5mm] {\scriptsize 1.10};
				\node	(v) at (11,4) [draw, circle, fill=lightred, inner sep=0pt,minimum size=5mm] {\scriptsize 0.80};
				\draw 
				(a) edge[->] (b) 
				(b) edge[->] (c)
				(c) edge[->] (d)
				(d) edge[->] (e)   
				(e) edge[->] node[near end,left] {\scriptsize 1.00} (f)
				(f) edge[->] node[near end,left] {\scriptsize 1.60} (g)
				(g) edge[->] node[near end,left] {\scriptsize 0.85} (h)
				(h) edge[->] (i)
				(i) edge[->] (j) 
				(j) edge[->] (k)
				(k) edge[->] node[near end,left] {\scriptsize 0.30} (l)
				(l) edge[->] node[near end,left] {\scriptsize 0.85} (m)   
				(m) edge[->] (n)
				(n) edge[->] (o)
				(o) edge[->] node[near end,left] {\scriptsize 0.30} (p)
				(p) edge[->] node[near end,left] {\scriptsize 0.85} (q)
				(q) edge[->] (r)
				(r) edge[->] (s)
				(s) edge[->] node[near end,left] {\scriptsize 0.30} (t)
				(t) edge[->] node[near end,left] {\scriptsize 0.85} (u)
				(u) edge[->] node[near end,left] {\scriptsize 1.35} (v)
				;
			\end{tikzpicture}
		}
		\caption{EP-2}
		\label{fig:ep-2}
	\end{subfigure}
	\caption{Computational structure of the evolved multigrid preconditioners. The color of the node denotes the type of operation. Black: Coarse-grid solver, Blue: Pointwise Jacobi smoothing, Red: Red-black Gauss-Seidel smoothing, White: No operation. The relaxation factor of each smoothing step is included in each node, while for coarse-grid correction, it is attached to the respective edge.}
	\label{fig:structure-evolved-preconditioners}
\end{figure}
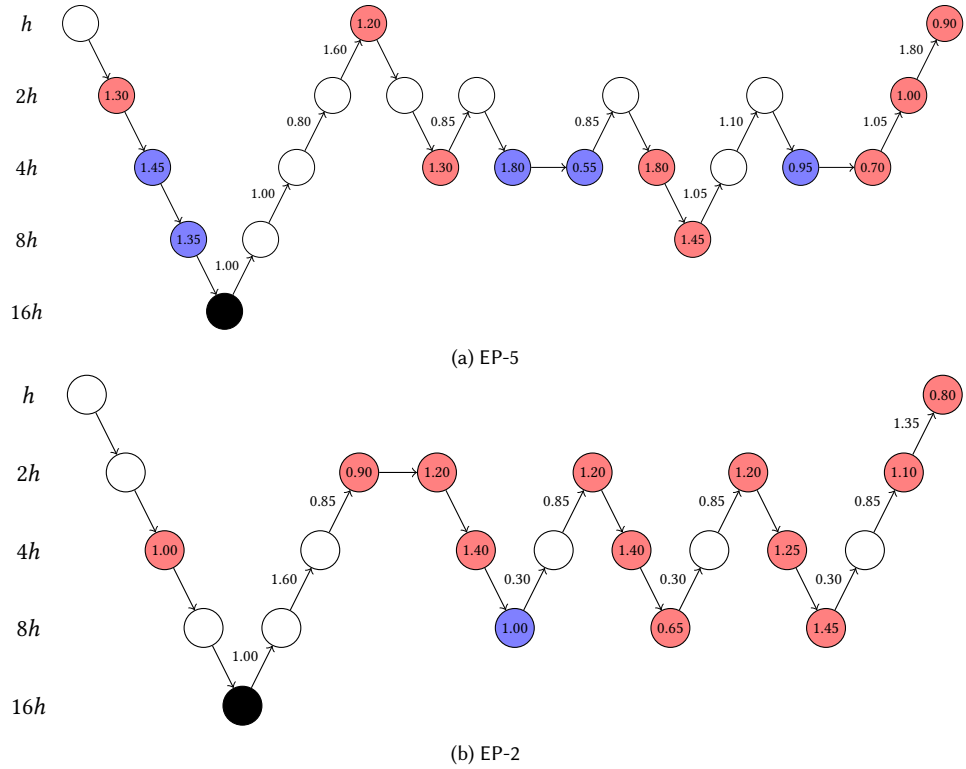
Figure~\ref{fig:structure-evolved-preconditioners} contains a graphical representation of each method's computational structure, including all relaxation factors.
These figures illustrate that, in each of the two cases, our grammatical representation of a multigrid preconditioner has enabled the construction of a unique sequence of computations, whose complexity exceeds those obtained by classical parameter optimization methods such as~\cite{oosterlee2003genetic,thekale2010optimizing,brown2021tuning}.
While both algorithms resemble a V-cycle, as the coarse-grid solver is only employed once, they include additional smoothing-based coarse-grid correction steps.
In contrast to classical multigrid cycles, these corrections are obtained from intermediate discretization levels without traversing the complete hierarchy down to the coarsest grid.
Furthermore, in both preconditioners, pre- and post-smoothing steps are omitted on certain levels, while the amount of smoothing is increased on others.
In EP-5 the number of smoothing steps is substantially higher on the coarser grids.
Since the computational cost of smoothing is significantly decreased with each coarsening step, this greatly reduces the overall execution time of the resulting preconditioned solver.
While EP-2 performs more smoothing on the second finest level, only a single step of red-black Gauss-Seidel is employed on the finest level.
In both preconditioners, red-black Gauss-Seidel is predominantly used for smoothing. 
However, especially EP-5 also includes intermediate Jacobi smoothing steps or a combination of both methods.
Finally, as shown in Figure~\ref{fig:structure-evolved-preconditioners}, both preconditioners combine a wide range of different relaxation factors within their computation.
\section{Conclusion and Future Work}
This work demonstrates how grammar-guided genetic programming (GGGP) can evolve multigrid preconditioners for Helmholtz problems that outperform known methods for different wavenumbers and even handle problems for which those methods fail.
Despite this accomplishment, further research is needed to investigate under which circumstances the presented approach can achieve consistent results.
We also aim to apply our approach to other multigrid variants, such as algebraic multigrid methods~\cite{xu2017algebraic}, and the solution of more challenging and complicated PDEs, such as nonlinear~\cite{henson2003multigrid} and saddle point problems~\cite{benzi2005numerical}.
Furthermore, the resulting implementation is mainly limited by the compute resources required to evaluate a sufficient number of preconditioners.
As a remedy, one could train a machine learning system to learn a model for predicting multigrid preconditioner performance based on the respective grammar representation.
Another promising research direction, which has been already mentioned in the introduction, is the grammar-based construction of a multigrid method from those components obtained by a machine learning-based optimization such as~\cite{greenfeld2019learning,huang2021learning}.
In addition, our approach could be utilized to accelerate the generation of training data for data-driven PDE solvers~\cite{li2020fourier,li2021physics,kovachki2021neural}. 
Finally, we aim to extend our implementation of GGGP to incorporate alternative initialization, crossover, and mutation operators, such as~\cite{garcia2006initialization,criado2020grammatically,couchet2007crossover}.
Also, while the implementation presented here employs tree-based GGGP, grammatical evolution (GE)~\cite{oneill2001grammatical,ryan2018handbook} represents a promising alternative, which could be as well integrated into our grammar-based approach for multigrid preconditioner design.
\clearpage
\bibliographystyle{ACM-Reference-Format}
\bibliography{references}


\begin{thebibliography}{55}


\ifx \showCODEN    \undefined \def \showCODEN     #1{\unskip}     \fi
\ifx \showDOI      \undefined \def \showDOI       #1{#1}\fi
\ifx \showISBNx    \undefined \def \showISBNx     #1{\unskip}     \fi
\ifx \showISBNxiii \undefined \def \showISBNxiii  #1{\unskip}     \fi
\ifx \showISSN     \undefined \def \showISSN      #1{\unskip}     \fi
\ifx \showLCCN     \undefined \def \showLCCN      #1{\unskip}     \fi
\ifx \shownote     \undefined \def \shownote      #1{#1}          \fi
\ifx \showarticletitle \undefined \def \showarticletitle #1{#1}   \fi
\ifx \showURL      \undefined \def \showURL       {\relax}        \fi
\providecommand\bibfield[2]{#2}
\providecommand\bibinfo[2]{#2}
\providecommand\natexlab[1]{#1}
\providecommand\showeprint[2][]{arXiv:#2}

\bibitem[\protect\citeauthoryear{Antoine and Darbas}{Antoine and
  Darbas}{2016}]%
        {antoine2016integral}
\bibfield{author}{\bibinfo{person}{Xavier Antoine} {and}
  \bibinfo{person}{Marion Darbas}.} \bibinfo{year}{2016}\natexlab{}.
\newblock \showarticletitle{{Integral Equations and Iterative Schemes for
  Acoustic Scattering Problems}}.
\newblock In \bibinfo{booktitle}{\emph{{Numerical Methods for Acoustics
  Problems}}}, \bibfield{editor}{\bibinfo{person}{F.~Magoul{\`e}s}} (Ed.).
  \bibinfo{publisher}{{Saxe-Coburg Editors}}.
\newblock
\urldef\tempurl%
\url{https://hal.archives-ouvertes.fr/hal-00591456}
\showURL{%
\tempurl}


\bibitem[\protect\citeauthoryear{Benzi}{Benzi}{2002}]%
        {benzi2002preconditioning}
\bibfield{author}{\bibinfo{person}{Michele Benzi}.}
  \bibinfo{year}{2002}\natexlab{}.
\newblock \showarticletitle{Preconditioning Techniques for Large Linear
  Systems: A Survey}.
\newblock \bibinfo{journal}{\emph{J. Comput. Phys.}} \bibinfo{volume}{182},
  \bibinfo{number}{2} (\bibinfo{year}{2002}), \bibinfo{pages}{418--477}.
\newblock
\showISSN{0021-9991}
\urldef\tempurl%
\url{https://doi.org/10.1006/jcph.2002.7176}
\showDOI{\tempurl}


\bibitem[\protect\citeauthoryear{Benzi, Golub, and Liesen}{Benzi
  et~al\mbox{.}}{2005}]%
        {benzi2005numerical}
\bibfield{author}{\bibinfo{person}{Michele Benzi}, \bibinfo{person}{Gene~H.
  Golub}, {and} \bibinfo{person}{Jörg Liesen}.}
  \bibinfo{year}{2005}\natexlab{}.
\newblock \showarticletitle{Numerical solution of saddle point problems}.
\newblock \bibinfo{journal}{\emph{Acta Numerica}}  \bibinfo{volume}{14}
  (\bibinfo{year}{2005}), \bibinfo{pages}{1–137}.
\newblock
\urldef\tempurl%
\url{https://doi.org/10.1017/S0962492904000212}
\showDOI{\tempurl}


\bibitem[\protect\citeauthoryear{Billette and Brandsberg-Dahl}{Billette and
  Brandsberg-Dahl}{2005}]%
        {billette20052004}
\bibfield{author}{\bibinfo{person}{F.~J. Billette} {and} \bibinfo{person}{S.
  Brandsberg-Dahl}.} \bibinfo{year}{2005}\natexlab{}.
\newblock \showarticletitle{The 2004 BP Velocity Benchmark}.
\newblock Article \bibinfo{articleno}{cp-1-00513}.
\newblock
\showISSN{2214-4609}
\urldef\tempurl%
\url{https://doi.org/10.3997/2214-4609-pdb.1.B035}
\showDOI{\tempurl}


\bibitem[\protect\citeauthoryear{Brandt}{Brandt}{1977}]%
        {brandt1977multi}
\bibfield{author}{\bibinfo{person}{Achi Brandt}.}
  \bibinfo{year}{1977}\natexlab{}.
\newblock \showarticletitle{Multi-level adaptive solutions to boundary-value
  problems}.
\newblock \bibinfo{journal}{\emph{Mathematics of computation}}
  \bibinfo{volume}{31}, \bibinfo{number}{138} (\bibinfo{year}{1977}),
  \bibinfo{pages}{333--390}.
\newblock


\bibitem[\protect\citeauthoryear{Briggs, Henson, and McCormick}{Briggs
  et~al\mbox{.}}{2000}]%
        {briggs2000multigrid}
\bibfield{author}{\bibinfo{person}{William~L. Briggs},
  \bibinfo{person}{Van~Emden Henson}, {and} \bibinfo{person}{Steve~F.
  McCormick}.} \bibinfo{year}{2000}\natexlab{}.
\newblock \bibinfo{booktitle}{\emph{A Multigrid Tutorial}
  (\bibinfo{edition}{second} ed.)}.
\newblock \bibinfo{publisher}{Society for Industrial and Applied Mathematics}.
\newblock
\urldef\tempurl%
\url{https://doi.org/10.1137/1.9780898719505}
\showDOI{\tempurl}


\bibitem[\protect\citeauthoryear{Brown, He, MacLachlan, Menickelly, and
  Wild}{Brown et~al\mbox{.}}{2021}]%
        {brown2021tuning}
\bibfield{author}{\bibinfo{person}{Jed Brown}, \bibinfo{person}{Yunhui He},
  \bibinfo{person}{Scott MacLachlan}, \bibinfo{person}{Matt Menickelly}, {and}
  \bibinfo{person}{Stefan~M. Wild}.} \bibinfo{year}{2021}\natexlab{}.
\newblock \showarticletitle{Tuning Multigrid Methods with Robust Optimization
  and Local Fourier Analysis}.
\newblock \bibinfo{journal}{\emph{SIAM Journal on Scientific Computing}}
  \bibinfo{volume}{43}, \bibinfo{number}{1} (\bibinfo{year}{2021}),
  \bibinfo{pages}{A109--A138}.
\newblock
\urldef\tempurl%
\url{https://doi.org/10.1137/19M1308669}
\showDOI{\tempurl}


\bibitem[\protect\citeauthoryear{Cocquet and Gander}{Cocquet and
  Gander}{2017}]%
        {cocquet2017shift}
\bibfield{author}{\bibinfo{person}{Pierre-Henri Cocquet} {and}
  \bibinfo{person}{Martin~J. Gander}.} \bibinfo{year}{2017}\natexlab{}.
\newblock \showarticletitle{How Large a Shift is Needed in the Shifted
  Helmholtz Preconditioner for its Effective Inversion by Multigrid?}
\newblock \bibinfo{journal}{\emph{SIAM Journal on Scientific Computing}}
  \bibinfo{volume}{39}, \bibinfo{number}{2} (\bibinfo{year}{2017}),
  \bibinfo{pages}{A438--A478}.
\newblock
\urldef\tempurl%
\url{https://doi.org/10.1137/15M102085X}
\showDOI{\tempurl}


\bibitem[\protect\citeauthoryear{Cools and Vanroose}{Cools and
  Vanroose}{2013}]%
        {cools2013analysis}
\bibfield{author}{\bibinfo{person}{Siegfried Cools} {and} \bibinfo{person}{Wim
  Vanroose}.} \bibinfo{year}{2013}\natexlab{}.
\newblock \showarticletitle{Local Fourier analysis of the complex shifted
  Laplacian preconditioner for Helmholtz problems}.
\newblock \bibinfo{journal}{\emph{Numerical Linear Algebra with Applications}}
  \bibinfo{volume}{20}, \bibinfo{number}{4} (\bibinfo{year}{2013}),
  \bibinfo{pages}{575--597}.
\newblock
\urldef\tempurl%
\url{https://doi.org/10.1002/nla.1881}
\showDOI{\tempurl}


\bibitem[\protect\citeauthoryear{Couchet, Manrique, R{\'i}os, and
  Rodr{\'i}guez-Pat{\'o}n}{Couchet et~al\mbox{.}}{2007}]%
        {couchet2007crossover}
\bibfield{author}{\bibinfo{person}{Jorge Couchet}, \bibinfo{person}{Daniel
  Manrique}, \bibinfo{person}{Juan R{\'i}os}, {and} \bibinfo{person}{Alfonso
  Rodr{\'i}guez-Pat{\'o}n}.} \bibinfo{year}{2007}\natexlab{}.
\newblock \showarticletitle{Crossover and mutation operators for grammar-guided
  genetic programming}.
\newblock \bibinfo{journal}{\emph{Soft Computing}} \bibinfo{volume}{11},
  \bibinfo{number}{10} (\bibinfo{date}{01 Aug} \bibinfo{year}{2007}),
  \bibinfo{pages}{943--955}.
\newblock
\showISSN{1433-7479}
\urldef\tempurl%
\url{https://doi.org/10.1007/s00500-006-0144-9}
\showDOI{\tempurl}


\bibitem[\protect\citeauthoryear{{Deb}, {Pratap}, {Agarwal}, and
  {Meyarivan}}{{Deb} et~al\mbox{.}}{2002}]%
        {deb2002fast}
\bibfield{author}{\bibinfo{person}{K. {Deb}}, \bibinfo{person}{A. {Pratap}},
  \bibinfo{person}{S. {Agarwal}}, {and} \bibinfo{person}{T. {Meyarivan}}.}
  \bibinfo{year}{2002}\natexlab{}.
\newblock \showarticletitle{A fast and elitist multiobjective genetic
  algorithm: {NSGA-II}}.
\newblock \bibinfo{journal}{\emph{IEEE Transactions on Evolutionary
  Computation}} \bibinfo{volume}{6}, \bibinfo{number}{2}
  (\bibinfo{year}{2002}), \bibinfo{pages}{182--197}.
\newblock
\urldef\tempurl%
\url{https://doi.org/10.1109/4235.996017}
\showDOI{\tempurl}


\bibitem[\protect\citeauthoryear{Erlangga}{Erlangga}{2008}]%
        {erlangga2008advances}
\bibfield{author}{\bibinfo{person}{Yogi~A. Erlangga}.}
  \bibinfo{year}{2008}\natexlab{}.
\newblock \showarticletitle{Advances in Iterative Methods and Preconditioners
  for the Helmholtz Equation}.
\newblock \bibinfo{journal}{\emph{Archives of Computational Methods in
  Engineering}} \bibinfo{volume}{15}, \bibinfo{number}{1} (\bibinfo{date}{01
  Mar} \bibinfo{year}{2008}), \bibinfo{pages}{37--66}.
\newblock
\showISSN{1886-1784}
\urldef\tempurl%
\url{https://doi.org/10.1007/s11831-007-9013-7}
\showDOI{\tempurl}


\bibitem[\protect\citeauthoryear{Erlangga, Oosterlee, and Vuik}{Erlangga
  et~al\mbox{.}}{2006}]%
        {erlangga2006multigrid}
\bibfield{author}{\bibinfo{person}{Y.~A. Erlangga}, \bibinfo{person}{C.~W.
  Oosterlee}, {and} \bibinfo{person}{C. Vuik}.}
  \bibinfo{year}{2006}\natexlab{}.
\newblock \showarticletitle{A Novel Multigrid Based Preconditioner For
  Heterogeneous Helmholtz Problems}.
\newblock \bibinfo{journal}{\emph{SIAM Journal on Scientific Computing}}
  \bibinfo{volume}{27}, \bibinfo{number}{4} (\bibinfo{year}{2006}),
  \bibinfo{pages}{1471--1492}.
\newblock
\urldef\tempurl%
\url{https://doi.org/10.1137/040615195}
\showDOI{\tempurl}


\bibitem[\protect\citeauthoryear{Erlangga, Vuik, and Oosterlee}{Erlangga
  et~al\mbox{.}}{2004}]%
        {erlangga2004preconditioner}
\bibfield{author}{\bibinfo{person}{Y.~A. Erlangga}, \bibinfo{person}{C. Vuik},
  {and} \bibinfo{person}{C.~W. Oosterlee}.} \bibinfo{year}{2004}\natexlab{}.
\newblock \showarticletitle{On a class of preconditioners for solving
  the Helmholtz equation}.
\newblock \bibinfo{journal}{\emph{Applied Numerical Mathematics}}
  \bibinfo{volume}{50}, \bibinfo{number}{3} (\bibinfo{year}{2004}),
  \bibinfo{pages}{409 -- 425}.
\newblock
\showISSN{0168-9274}
\urldef\tempurl%
\url{https://doi.org/10.1016/j.apnum.2004.01.009}
\showDOI{\tempurl}


\bibitem[\protect\citeauthoryear{Ernst and Gander}{Ernst and Gander}{2012}]%
        {ernst2012difficult}
\bibfield{author}{\bibinfo{person}{O.~G. Ernst} {and} \bibinfo{person}{M.~J.
  Gander}.} \bibinfo{year}{2012}\natexlab{}.
\newblock \bibinfo{booktitle}{\emph{Why it is Difficult to Solve Helmholtz
  Problems with Classical Iterative Methods}}.
\newblock \bibinfo{publisher}{Springer Berlin Heidelberg},
  \bibinfo{address}{Berlin, Heidelberg}, \bibinfo{pages}{325--363}.
\newblock
\showISBNx{978-3-642-22061-6}
\urldef\tempurl%
\url{https://doi.org/10.1007/978-3-642-22061-6_10}
\showDOI{\tempurl}


\bibitem[\protect\citeauthoryear{Fedorenko}{Fedorenko}{1962}]%
        {fedorenko1962relaxation}
\bibfield{author}{\bibinfo{person}{Radii~Petrovich Fedorenko}.}
  \bibinfo{year}{1962}\natexlab{}.
\newblock \showarticletitle{A relaxation method for solving elliptic difference
  equations}.
\newblock \bibinfo{journal}{\emph{U. S. S. R. Comput. Math. and Math. Phys.}}
  \bibinfo{volume}{1}, \bibinfo{number}{4} (\bibinfo{year}{1962}),
  \bibinfo{pages}{1092--1096}.
\newblock
\showISSN{0041-5553}
\urldef\tempurl%
\url{https://doi.org/10.1016/0041-5553(62)90031-9}
\showDOI{\tempurl}


\bibitem[\protect\citeauthoryear{Fortin, Grenier, and Parizeau}{Fortin
  et~al\mbox{.}}{2013}]%
        {fortin2013generalizing}
\bibfield{author}{\bibinfo{person}{F\'{e}lix-Antoine Fortin},
  \bibinfo{person}{Simon Grenier}, {and} \bibinfo{person}{Marc Parizeau}.}
  \bibinfo{year}{2013}\natexlab{}.
\newblock \showarticletitle{Generalizing the Improved Run-Time Complexity
  Algorithm for Non-Dominated Sorting}. In
  \bibinfo{booktitle}{\emph{Proceedings of the 15th Annual Conference on
  Genetic and Evolutionary Computation}} (Amsterdam, The Netherlands)
  \emph{(\bibinfo{series}{GECCO '13})}. \bibinfo{publisher}{Association for
  Computing Machinery}, \bibinfo{address}{New York, NY, USA},
  \bibinfo{pages}{615–622}.
\newblock
\showISBNx{9781450319638}
\urldef\tempurl%
\url{https://doi.org/10.1145/2463372.2463454}
\showDOI{\tempurl}


\bibitem[\protect\citeauthoryear{Gander and Zhang}{Gander and Zhang}{2019}]%
        {gander2019solvers}
\bibfield{author}{\bibinfo{person}{Martin~J. Gander} {and} \bibinfo{person}{Hui
  Zhang}.} \bibinfo{year}{2019}\natexlab{}.
\newblock \showarticletitle{A Class of Iterative Solvers for the Helmholtz
  Equation: Factorizations, Sweeping Preconditioners, Source Transfer, Single
  Layer Potentials, Polarized Traces, and Optimized Schwarz Methods}.
\newblock \bibinfo{journal}{\emph{SIAM Rev.}} \bibinfo{volume}{61},
  \bibinfo{number}{1} (\bibinfo{year}{2019}), \bibinfo{pages}{3--76}.
\newblock
\urldef\tempurl%
\url{https://doi.org/10.1137/16M109781X}
\showDOI{\tempurl}


\bibitem[\protect\citeauthoryear{Garc{\'i}a-Arnau, Manrique, R{\'i}os, and
  Rodr{\'i}guez-Pat{\'o}n}{Garc{\'i}a-Arnau et~al\mbox{.}}{2007}]%
        {garcia2006initialization}
\bibfield{author}{\bibinfo{person}{M. Garc{\'i}a-Arnau}, \bibinfo{person}{D.
  Manrique}, \bibinfo{person}{J. R{\'i}os}, {and} \bibinfo{person}{A.
  Rodr{\'i}guez-Pat{\'o}n}.} \bibinfo{year}{2007}\natexlab{}.
\newblock \showarticletitle{Initialization Method for Grammar-Guided Genetic
  Programming}. In \bibinfo{booktitle}{\emph{Research and Development in
  Intelligent Systems XXIII}}, \bibfield{editor}{\bibinfo{person}{Max Bramer},
  \bibinfo{person}{Frans Coenen}, {and} \bibinfo{person}{Andrew Tuson}} (Eds.).
  \bibinfo{publisher}{Springer London}, \bibinfo{address}{London},
  \bibinfo{pages}{32--44}.
\newblock
\showISBNx{978-1-84628-663-6}
\urldef\tempurl%
\url{https://doi.org/10.1007/978-1-84628-663-6_3}
\showDOI{\tempurl}


\bibitem[\protect\citeauthoryear{Greenfeld, Galun, Basri, Yavneh, and
  Kimmel}{Greenfeld et~al\mbox{.}}{2019}]%
        {greenfeld2019learning}
\bibfield{author}{\bibinfo{person}{Daniel Greenfeld}, \bibinfo{person}{Meirav
  Galun}, \bibinfo{person}{Ronen Basri}, \bibinfo{person}{Irad Yavneh}, {and}
  \bibinfo{person}{Ron Kimmel}.} \bibinfo{year}{2019}\natexlab{}.
\newblock \showarticletitle{Learning to Optimize Multigrid {PDE} Solvers}. In
  \bibinfo{booktitle}{\emph{Proceedings of the 36th International Conference on
  Machine Learning}} \emph{(\bibinfo{series}{Proceedings of Machine Learning
  Research}, Vol.~\bibinfo{volume}{97})},
  \bibfield{editor}{\bibinfo{person}{Kamalika Chaudhuri} {and}
  \bibinfo{person}{Ruslan Salakhutdinov}} (Eds.). \bibinfo{publisher}{PMLR},
  \bibinfo{pages}{2415--2423}.
\newblock
\urldef\tempurl%
\url{https://proceedings.mlr.press/v97/greenfeld19a.html}
\showURL{%
\tempurl}


\bibitem[\protect\citeauthoryear{Hackbusch}{Hackbusch}{1985}]%
        {hackbusch2013multi}
\bibfield{author}{\bibinfo{person}{Wolfgang Hackbusch}.}
  \bibinfo{year}{1985}\natexlab{}.
\newblock \bibinfo{booktitle}{\emph{Multi-Grid Methods and Applications}}.
\newblock \bibinfo{publisher}{Springer-Verlag}.
\newblock


\bibitem[\protect\citeauthoryear{Henson}{Henson}{2003}]%
        {henson2003multigrid}
\bibfield{author}{\bibinfo{person}{Van~Emden Henson}.}
  \bibinfo{year}{2003}\natexlab{}.
\newblock \showarticletitle{{Multigrid methods nonlinear problems: an
  overview}}. In \bibinfo{booktitle}{\emph{Computational Imaging}},
  \bibfield{editor}{\bibinfo{person}{Charles~A. Bouman} {and}
  \bibinfo{person}{Robert~L. Stevenson}} (Eds.), Vol.~\bibinfo{volume}{5016}.
  International Society for Optics and Photonics, \bibinfo{publisher}{SPIE},
  \bibinfo{pages}{36 -- 48}.
\newblock
\urldef\tempurl%
\url{https://doi.org/10.1117/12.499473}
\showDOI{\tempurl}


\bibitem[\protect\citeauthoryear{Huang, Li, and Xi}{Huang
  et~al\mbox{.}}{2021}]%
        {huang2021learning}
\bibfield{author}{\bibinfo{person}{Ru Huang}, \bibinfo{person}{Ruipeng Li},
  {and} \bibinfo{person}{Yuanzhe Xi}.} \bibinfo{year}{2021}\natexlab{}.
\newblock \showarticletitle{Learning optimal multigrid smoothers via neural
  networks}.
\newblock  (\bibinfo{year}{2021}).
\newblock
\urldef\tempurl%
\url{https://doi.org/10.48550/ARXIV.2102.12071}
\showDOI{\tempurl}


\bibitem[\protect\citeauthoryear{Karniadakis, Kevrekidis, Lu, Perdikaris, Wang,
  and Yang}{Karniadakis et~al\mbox{.}}{2021}]%
        {karniadakis2021physics}
\bibfield{author}{\bibinfo{person}{George~Em Karniadakis},
  \bibinfo{person}{Ioannis~G. Kevrekidis}, \bibinfo{person}{Lu Lu},
  \bibinfo{person}{Paris Perdikaris}, \bibinfo{person}{Sifan Wang}, {and}
  \bibinfo{person}{Liu Yang}.} \bibinfo{year}{2021}\natexlab{}.
\newblock \showarticletitle{Physics-informed machine learning}.
\newblock \bibinfo{journal}{\emph{Nature Reviews Physics}} \bibinfo{volume}{3},
  \bibinfo{number}{6} (\bibinfo{date}{01 Jun} \bibinfo{year}{2021}),
  \bibinfo{pages}{422--440}.
\newblock
\showISSN{2522-5820}
\urldef\tempurl%
\url{https://doi.org/10.1038/s42254-021-00314-5}
\showDOI{\tempurl}


\bibitem[\protect\citeauthoryear{Katrutsa, Daulbaev, and Oseledets}{Katrutsa
  et~al\mbox{.}}{2020}]%
        {katrutsa2020black}
\bibfield{author}{\bibinfo{person}{Alexandr Katrutsa}, \bibinfo{person}{Talgat
  Daulbaev}, {and} \bibinfo{person}{Ivan Oseledets}.}
  \bibinfo{year}{2020}\natexlab{}.
\newblock \showarticletitle{Black-box learning of multigrid parameters}.
\newblock \bibinfo{journal}{\emph{J. Comput. Appl. Math.}}
  \bibinfo{volume}{368} (\bibinfo{year}{2020}), \bibinfo{pages}{112524}.
\newblock
\showISSN{0377-0427}
\urldef\tempurl%
\url{https://doi.org/10.1016/j.cam.2019.112524}
\showDOI{\tempurl}


\bibitem[\protect\citeauthoryear{K{\"o}stler, Heisig, Kohl, Kuckuk, Bauer, and
  R{\"u}de}{K{\"o}stler et~al\mbox{.}}{2020}]%
        {kostler2020code}
\bibfield{author}{\bibinfo{person}{Harald K{\"o}stler}, \bibinfo{person}{Marco
  Heisig}, \bibinfo{person}{Nils Kohl}, \bibinfo{person}{Sebastian Kuckuk},
  \bibinfo{person}{Martin Bauer}, {and} \bibinfo{person}{Ulrich R{\"u}de}.}
  \bibinfo{year}{2020}\natexlab{}.
\newblock \showarticletitle{Code generation approaches for parallel geometric
  multigrid solvers}.
\newblock \bibinfo{journal}{\emph{Analele Universitatii Ovidius Constanta -
  Seria Matematica}} \bibinfo{volume}{28}, \bibinfo{number}{3}
  (\bibinfo{year}{2020}), \bibinfo{pages}{123--152}.
\newblock


\bibitem[\protect\citeauthoryear{Kovachki, Li, Liu, Azizzadenesheli,
  Bhattacharya, Stuart, and Anandkumar}{Kovachki et~al\mbox{.}}{2021}]%
        {kovachki2021neural}
\bibfield{author}{\bibinfo{person}{Nikola Kovachki}, \bibinfo{person}{Zongyi
  Li}, \bibinfo{person}{Burigede Liu}, \bibinfo{person}{Kamyar
  Azizzadenesheli}, \bibinfo{person}{Kaushik Bhattacharya},
  \bibinfo{person}{Andrew Stuart}, {and} \bibinfo{person}{Anima Anandkumar}.}
  \bibinfo{year}{2021}\natexlab{}.
\newblock \showarticletitle{Neural Operator: Learning Maps Between Function
  Spaces}.
\newblock  (\bibinfo{year}{2021}).
\newblock
\urldef\tempurl%
\url{https://doi.org/10.48550/ARXIV.2108.08481}
\showDOI{\tempurl}


\bibitem[\protect\citeauthoryear{Koza}{Koza}{1994}]%
        {koza1992genetic}
\bibfield{author}{\bibinfo{person}{John~R. Koza}.}
  \bibinfo{year}{1994}\natexlab{}.
\newblock \bibinfo{booktitle}{\emph{Genetic programming as a means for
  programming computers by natural selection}}. Vol.~\bibinfo{volume}{4}.
\newblock 87--112 pages.
\newblock
\showISSN{1573-1375}
\urldef\tempurl%
\url{https://doi.org/10.1007/BF00175355}
\showDOI{\tempurl}


\bibitem[\protect\citeauthoryear{Köstler and Rüde}{Köstler and
  Rüde}{2004}]%
        {koestler2004extrapolation}
\bibfield{author}{\bibinfo{person}{H. Köstler} {and} \bibinfo{person}{U.
  Rüde}.} \bibinfo{year}{2004}\natexlab{}.
\newblock \showarticletitle{Extrapolation Techniques for Computing Accurate
  Solutions of Elliptic Problems with Singular Solutions}. In
  \bibinfo{booktitle}{\emph{Computational Science - ICCS 2004}},
  \bibfield{editor}{\bibinfo{person}{Marian Bubak}, \bibinfo{person}{Geert~Dick
  van Albada}, \bibinfo{person}{Peter M.~A. Sloot}, {and} \bibinfo{person}{Jack
  Dongarra}} (Eds.). \bibinfo{publisher}{Springer Berlin Heidelberg},
  \bibinfo{address}{Berlin, Heidelberg}, \bibinfo{pages}{410--417}.
\newblock
\showISBNx{978-3-540-25944-2}
\urldef\tempurl%
\url{https://doi.org/10.1007/978-3-540-25944-2_54}
\showURL{%
\tempurl}


\bibitem[\protect\citeauthoryear{Lengauer, Apel, Bolten, Chiba, R{\"u}de,
  Teich, Gr{\"o}{\ss}linger, Hannig, K{\"o}stler, Claus,
  et~al\mbox{.}}{Lengauer et~al\mbox{.}}{2020}]%
        {lengauer2020exastencils}
\bibfield{author}{\bibinfo{person}{Christian Lengauer}, \bibinfo{person}{Sven
  Apel}, \bibinfo{person}{Matthias Bolten}, \bibinfo{person}{Shigeru Chiba},
  \bibinfo{person}{Ulrich R{\"u}de}, \bibinfo{person}{J{\"u}rgen Teich},
  \bibinfo{person}{Armin Gr{\"o}{\ss}linger}, \bibinfo{person}{Frank Hannig},
  \bibinfo{person}{Harald K{\"o}stler}, \bibinfo{person}{Lisa Claus},
  {et~al\mbox{.}}} \bibinfo{year}{2020}\natexlab{}.
\newblock \showarticletitle{ExaStencils: Advanced Multigrid Solver Generation}.
  In \bibinfo{booktitle}{\emph{Software for Exascale Computing - SPPEXA
  2016-2019}}, \bibfield{editor}{\bibinfo{person}{Hans-Joachim Bungartz},
  \bibinfo{person}{Severin Reiz}, \bibinfo{person}{Benjamin Uekermann},
  \bibinfo{person}{Philipp Neumann}, {and} \bibinfo{person}{Wolfgang~E. Nagel}}
  (Eds.). \bibinfo{publisher}{Springer International Publishing},
  \bibinfo{address}{Cham}, \bibinfo{pages}{405--452}.
\newblock
\showISBNx{978-3-030-47956-5}
\urldef\tempurl%
\url{https://doi.org/10.1007/978-3-030-47956-5_14}
\showURL{%
\tempurl}


\bibitem[\protect\citeauthoryear{Li, Kovachki, Azizzadenesheli, Liu,
  Bhattacharya, Stuart, and Anandkumar}{Li et~al\mbox{.}}{2020}]%
        {li2020fourier}
\bibfield{author}{\bibinfo{person}{Zongyi Li}, \bibinfo{person}{Nikola
  Kovachki}, \bibinfo{person}{Kamyar Azizzadenesheli},
  \bibinfo{person}{Burigede Liu}, \bibinfo{person}{Kaushik Bhattacharya},
  \bibinfo{person}{Andrew Stuart}, {and} \bibinfo{person}{Anima Anandkumar}.}
  \bibinfo{year}{2020}\natexlab{}.
\newblock \showarticletitle{Fourier Neural Operator for Parametric Partial
  Differential Equations}.
\newblock  (\bibinfo{year}{2020}).
\newblock
\urldef\tempurl%
\url{https://doi.org/10.48550/ARXIV.2010.08895}
\showDOI{\tempurl}


\bibitem[\protect\citeauthoryear{Li, Zheng, Kovachki, Jin, Chen, Liu,
  Azizzadenesheli, and Anandkumar}{Li et~al\mbox{.}}{2021}]%
        {li2021physics}
\bibfield{author}{\bibinfo{person}{Zongyi Li}, \bibinfo{person}{Hongkai Zheng},
  \bibinfo{person}{Nikola Kovachki}, \bibinfo{person}{David Jin},
  \bibinfo{person}{Haoxuan Chen}, \bibinfo{person}{Burigede Liu},
  \bibinfo{person}{Kamyar Azizzadenesheli}, {and} \bibinfo{person}{Anima
  Anandkumar}.} \bibinfo{year}{2021}\natexlab{}.
\newblock \showarticletitle{Physics-Informed Neural Operator for Learning
  Partial Differential Equations}.
\newblock  (\bibinfo{year}{2021}).
\newblock
\urldef\tempurl%
\url{https://doi.org/10.48550/ARXIV.2111.03794}
\showDOI{\tempurl}


\bibitem[\protect\citeauthoryear{Linz and Rodger}{Linz and Rodger}{2022}]%
        {linz2006introduction}
\bibfield{author}{\bibinfo{person}{Peter Linz} {and} \bibinfo{person}{Susan~H.
  Rodger}.} \bibinfo{year}{2022}\natexlab{}.
\newblock \bibinfo{booktitle}{\emph{An Introduction to Formal Languages and
  Automata}}.
\newblock \bibinfo{publisher}{Jones \& Bartlett Learning}.
\newblock


\bibitem[\protect\citeauthoryear{Luz, Galun, Maron, Basri, and Yavneh}{Luz
  et~al\mbox{.}}{2020}]%
        {luz2020learning}
\bibfield{author}{\bibinfo{person}{Ilay Luz}, \bibinfo{person}{Meirav Galun},
  \bibinfo{person}{Haggai Maron}, \bibinfo{person}{Ronen Basri}, {and}
  \bibinfo{person}{Irad Yavneh}.} \bibinfo{year}{2020}\natexlab{}.
\newblock \showarticletitle{Learning Algebraic Multigrid Using Graph Neural
  Networks}. In \bibinfo{booktitle}{\emph{Proceedings of the 37th International
  Conference on Machine Learning}} \emph{(\bibinfo{series}{Proceedings of
  Machine Learning Research}, Vol.~\bibinfo{volume}{119})},
  \bibfield{editor}{\bibinfo{person}{Hal~Daumé III} {and}
  \bibinfo{person}{Aarti Singh}} (Eds.). \bibinfo{publisher}{PMLR},
  \bibinfo{pages}{6489--6499}.
\newblock
\urldef\tempurl%
\url{https://proceedings.mlr.press/v119/luz20a.html}
\showURL{%
\tempurl}


\bibitem[\protect\citeauthoryear{Martin, Wiley, and Marfurt}{Martin
  et~al\mbox{.}}{2006}]%
        {martin2006marmousi2}
\bibfield{author}{\bibinfo{person}{Gary~S. Martin}, \bibinfo{person}{Robert
  Wiley}, {and} \bibinfo{person}{Kurt~J. Marfurt}.}
  \bibinfo{year}{2006}\natexlab{}.
\newblock \showarticletitle{Marmousi2: An elastic upgrade for Marmousi}.
\newblock \bibinfo{journal}{\emph{The Leading Edge}} \bibinfo{volume}{25},
  \bibinfo{number}{2} (\bibinfo{year}{2006}), \bibinfo{pages}{156--166}.
\newblock
\urldef\tempurl%
\url{https://doi.org/10.1190/1.2172306}
\showDOI{\tempurl}


\bibitem[\protect\citeauthoryear{McKay, Hoai, Whigham, Shan, and O'Neill}{McKay
  et~al\mbox{.}}{2010}]%
        {mckay2010grammar}
\bibfield{author}{\bibinfo{person}{Robert~I. McKay},
  \bibinfo{person}{Nguyen~Xuan Hoai}, \bibinfo{person}{Peter~Alexander
  Whigham}, \bibinfo{person}{Yin Shan}, {and} \bibinfo{person}{Michael
  O'Neill}.} \bibinfo{year}{2010}\natexlab{}.
\newblock \showarticletitle{Grammar-based Genetic Programming: a survey}.
\newblock \bibinfo{journal}{\emph{Genetic Programming and Evolvable Machines}}
  \bibinfo{volume}{11}, \bibinfo{number}{3} (\bibinfo{date}{01 Sep}
  \bibinfo{year}{2010}), \bibinfo{pages}{365--396}.
\newblock
\showISSN{1573-7632}
\urldef\tempurl%
\url{https://doi.org/10.1007/s10710-010-9109-y}
\showDOI{\tempurl}


\bibitem[\protect\citeauthoryear{O'Neill and Ryan}{O'Neill and Ryan}{2001}]%
        {oneill2001grammatical}
\bibfield{author}{\bibinfo{person}{M. O'Neill} {and} \bibinfo{person}{C.
  Ryan}.} \bibinfo{year}{2001}\natexlab{}.
\newblock \showarticletitle{Grammatical evolution}.
\newblock \bibinfo{journal}{\emph{IEEE Transactions on Evolutionary
  Computation}} \bibinfo{volume}{5}, \bibinfo{number}{4}
  (\bibinfo{year}{2001}), \bibinfo{pages}{349--358}.
\newblock
\urldef\tempurl%
\url{https://doi.org/10.1109/4235.942529}
\showDOI{\tempurl}


\bibitem[\protect\citeauthoryear{Oosterlee and Wienands}{Oosterlee and
  Wienands}{2003}]%
        {oosterlee2003genetic}
\bibfield{author}{\bibinfo{person}{C.~W. Oosterlee} {and} \bibinfo{person}{R.
  Wienands}.} \bibinfo{year}{2003}\natexlab{}.
\newblock \showarticletitle{A Genetic Search for Optimal Multigrid Components
  Within a Fourier Analysis Setting}.
\newblock \bibinfo{journal}{\emph{SIAM Journal on Scientific Computing}}
  \bibinfo{volume}{24}, \bibinfo{number}{3} (\bibinfo{year}{2003}),
  \bibinfo{pages}{924--944}.
\newblock
\urldef\tempurl%
\url{https://doi.org/10.1137/S1064827501397950}
\showDOI{\tempurl}


\bibitem[\protect\citeauthoryear{Osei-Kuffuor and Saad}{Osei-Kuffuor and
  Saad}{2010}]%
        {oseikuffuor2010preconditioning}
\bibfield{author}{\bibinfo{person}{Daniel Osei-Kuffuor} {and}
  \bibinfo{person}{Yousef Saad}.} \bibinfo{year}{2010}\natexlab{}.
\newblock \showarticletitle{Preconditioning Helmholtz linear systems}.
\newblock \bibinfo{journal}{\emph{Applied Numerical Mathematics}}
  \bibinfo{volume}{60}, \bibinfo{number}{4} (\bibinfo{year}{2010}),
  \bibinfo{pages}{420 -- 431}.
\newblock
\showISSN{0168-9274}
\urldef\tempurl%
\url{https://doi.org/10.1016/j.apnum.2009.09.003}
\showDOI{\tempurl}
\newblock
\shownote{Special Issue: NUMAN 2008}.


\bibitem[\protect\citeauthoryear{Pierce}{Pierce}{2019}]%
        {pierce2019acoustics}
\bibfield{author}{\bibinfo{person}{Allan~D. Pierce}.}
  \bibinfo{year}{2019}\natexlab{}.
\newblock \bibinfo{booktitle}{\emph{Acoustics: An Introduction to Its Physical
  Principles and Applications}}.
\newblock \bibinfo{publisher}{Springer}.
\newblock
\urldef\tempurl%
\url{10.1007/978-3-030-11214-1}
\showURL{%
\tempurl}


\bibitem[\protect\citeauthoryear{Poli, Langdon, and McPhee}{Poli
  et~al\mbox{.}}{2008}]%
        {poli2008field}
\bibfield{author}{\bibinfo{person}{Riccardo Poli}, \bibinfo{person}{William~B.
  Langdon}, {and} \bibinfo{person}{Nicholas~Freitag McPhee}.}
  \bibinfo{year}{2008}\natexlab{}.
\newblock \bibinfo{booktitle}{\emph{A Field Guide to Genetic Programming}}.
\newblock \bibinfo{publisher}{Lulu Enterprises, UK Ltd}.
\newblock
\showISBNx{1409200736}


\bibitem[\protect\citeauthoryear{Raissi, Perdikaris, and Karniadakis}{Raissi
  et~al\mbox{.}}{2019}]%
        {raissi2019physics}
\bibfield{author}{\bibinfo{person}{M. Raissi}, \bibinfo{person}{P. Perdikaris},
  {and} \bibinfo{person}{G.~E. Karniadakis}.} \bibinfo{year}{2019}\natexlab{}.
\newblock \showarticletitle{Physics-informed neural networks: A deep learning
  framework for solving forward and inverse problems involving nonlinear
  partial differential equations}.
\newblock \bibinfo{journal}{\emph{J. Comput. Phys.}}  \bibinfo{volume}{378}
  (\bibinfo{year}{2019}), \bibinfo{pages}{686--707}.
\newblock
\showISSN{0021-9991}
\urldef\tempurl%
\url{https://doi.org/10.1016/j.jcp.2018.10.045}
\showDOI{\tempurl}


\bibitem[\protect\citeauthoryear{Ramos~Criado, Barrios~Rolan{\'i}a, Manrique,
  and Serrano}{Ramos~Criado et~al\mbox{.}}{2020}]%
        {criado2020grammatically}
\bibfield{author}{\bibinfo{person}{Pablo Ramos~Criado}, \bibinfo{person}{D.
  Barrios~Rolan{\'i}a}, \bibinfo{person}{Daniel Manrique}, {and}
  \bibinfo{person}{Emilio Serrano}.} \bibinfo{year}{2020}\natexlab{}.
\newblock \showarticletitle{Grammatically uniform population initialization for
  grammar-guided genetic programming}.
\newblock \bibinfo{journal}{\emph{Soft Computing}} \bibinfo{volume}{24},
  \bibinfo{number}{15} (\bibinfo{date}{01 Aug} \bibinfo{year}{2020}),
  \bibinfo{pages}{11265--11282}.
\newblock
\showISSN{1433-7479}
\urldef\tempurl%
\url{https://doi.org/10.1007/s00500-020-05061-w}
\showDOI{\tempurl}


\bibitem[\protect\citeauthoryear{Ryan, O'Neill, and Collins}{Ryan
  et~al\mbox{.}}{2018}]%
        {ryan2018handbook}
\bibfield{author}{\bibinfo{person}{Conor Ryan}, \bibinfo{person}{Michael
  O'Neill}, {and} \bibinfo{person}{J.~J. Collins}.}
  \bibinfo{year}{2018}\natexlab{}.
\newblock \bibinfo{booktitle}{\emph{Handbook of Grammatical Evolution}}.
\newblock \bibinfo{publisher}{Springer}.
\newblock
\urldef\tempurl%
\url{https://doi.org/10.1007/978-3-319-78717-6}
\showDOI{\tempurl}


\bibitem[\protect\citeauthoryear{Saad}{Saad}{2003}]%
        {saad2003iterative}
\bibfield{author}{\bibinfo{person}{Yousef Saad}.}
  \bibinfo{year}{2003}\natexlab{}.
\newblock \bibinfo{booktitle}{\emph{Iterative Methods for Sparse Linear
  Systems} (\bibinfo{edition}{second} ed.)}.
\newblock \bibinfo{publisher}{Society for Industrial and Applied Mathematics}.
\newblock
\urldef\tempurl%
\url{https://doi.org/10.1137/1.9780898718003}
\showDOI{\tempurl}


\bibitem[\protect\citeauthoryear{Salon and Chari}{Salon and Chari}{1999}]%
        {salon1999numerical}
\bibfield{author}{\bibinfo{person}{Sheppard Salon} {and} \bibinfo{person}{M.
  Chari}.} \bibinfo{year}{1999}\natexlab{}.
\newblock \bibinfo{booktitle}{\emph{Numerical methods in electromagnetism}}.
\newblock \bibinfo{publisher}{Elsevier}.
\newblock


\bibitem[\protect\citeauthoryear{Schmitt, Kuckuk, Hannig, Köstler, and
  Teich}{Schmitt et~al\mbox{.}}{2014}]%
        {schmitt2014exaslang}
\bibfield{author}{\bibinfo{person}{Christian Schmitt},
  \bibinfo{person}{Sebastian Kuckuk}, \bibinfo{person}{Frank Hannig},
  \bibinfo{person}{Harald Köstler}, {and} \bibinfo{person}{Jürgen Teich}.}
  \bibinfo{year}{2014}\natexlab{}.
\newblock \showarticletitle{ExaSlang: A Domain-Specific Language for Highly
  Scalable Multigrid Solvers}. In \bibinfo{booktitle}{\emph{2014 Fourth
  International Workshop on Domain-Specific Languages and High-Level Frameworks
  for High Performance Computing}}. \bibinfo{pages}{42--51}.
\newblock
\urldef\tempurl%
\url{https://doi.org/10.1109/WOLFHPC.2014.11}
\showDOI{\tempurl}


\bibitem[\protect\citeauthoryear{Schmitt, Kuckuk, and K\"{o}stler}{Schmitt
  et~al\mbox{.}}{2020}]%
        {schmitt2020constructing}
\bibfield{author}{\bibinfo{person}{Jonas Schmitt}, \bibinfo{person}{Sebastian
  Kuckuk}, {and} \bibinfo{person}{Harald K\"{o}stler}.}
  \bibinfo{year}{2020}\natexlab{}.
\newblock \showarticletitle{Constructing Efficient Multigrid Solvers with
  Genetic Programming}. In \bibinfo{booktitle}{\emph{Proceedings of the 2020
  Genetic and Evolutionary Computation Conference}} (Canc\'{u}n, Mexico)
  \emph{(\bibinfo{series}{GECCO '20})}. \bibinfo{publisher}{Association for
  Computing Machinery}, \bibinfo{address}{New York, NY, USA},
  \bibinfo{pages}{1012–1020}.
\newblock
\showISBNx{9781450371285}
\urldef\tempurl%
\url{https://doi.org/10.1145/3377930.3389811}
\showDOI{\tempurl}


\bibitem[\protect\citeauthoryear{Schmitt, Kuckuk, and K{\"o}stler}{Schmitt
  et~al\mbox{.}}{2021}]%
        {schmitt2021evostencils}
\bibfield{author}{\bibinfo{person}{Jonas Schmitt}, \bibinfo{person}{Sebastian
  Kuckuk}, {and} \bibinfo{person}{Harald K{\"o}stler}.}
  \bibinfo{year}{2021}\natexlab{}.
\newblock \showarticletitle{EvoStencils: a grammar-based genetic programming
  approach for constructing efficient geometric multigrid methods}.
\newblock \bibinfo{journal}{\emph{Genetic Programming and Evolvable Machines}}
  (\bibinfo{date}{03 Sep} \bibinfo{year}{2021}).
\newblock
\showISSN{1573-7632}
\urldef\tempurl%
\url{https://doi.org/10.1007/s10710-021-09412-w}
\showDOI{\tempurl}


\bibitem[\protect\citeauthoryear{Taghibakhshi, MacLachlan, Olson, and
  West}{Taghibakhshi et~al\mbox{.}}{2021}]%
        {taghibakhshi2021optimization}
\bibfield{author}{\bibinfo{person}{Ali Taghibakhshi}, \bibinfo{person}{Scott
  MacLachlan}, \bibinfo{person}{Luke Olson}, {and} \bibinfo{person}{Matthew
  West}.} \bibinfo{year}{2021}\natexlab{}.
\newblock \showarticletitle{Optimization-Based Algebraic Multigrid Coarsening
  Using Reinforcement Learning}.
\newblock   \bibinfo{volume}{34} (\bibinfo{year}{2021}),
  \bibinfo{pages}{12129--12140}.
\newblock
\urldef\tempurl%
\url{https://proceedings.neurips.cc/paper/2021/file/6531b32f8d02fece98ff36a64a7c8260-Paper.pdf}
\showURL{%
\tempurl}


\bibitem[\protect\citeauthoryear{Thekale, Gradl, Klamroth, and Rüde}{Thekale
  et~al\mbox{.}}{2010}]%
        {thekale2010optimizing}
\bibfield{author}{\bibinfo{person}{A. Thekale}, \bibinfo{person}{T. Gradl},
  \bibinfo{person}{K. Klamroth}, {and} \bibinfo{person}{U. Rüde}.}
  \bibinfo{year}{2010}\natexlab{}.
\newblock \showarticletitle{Optimizing the number of multigrid cycles in the
  full multigrid algorithm}.
\newblock \bibinfo{journal}{\emph{Numerical Linear Algebra with Applications}}
  \bibinfo{volume}{17}, \bibinfo{number}{2‐3} (\bibinfo{year}{2010}),
  \bibinfo{pages}{199--210}.
\newblock
\urldef\tempurl%
\url{https://doi.org/10.1002/nla.697}
\showDOI{\tempurl}


\bibitem[\protect\citeauthoryear{Trottenberg, Oosterlee, and
  Schuller}{Trottenberg et~al\mbox{.}}{2000}]%
        {trottenberg2000multigrid}
\bibfield{author}{\bibinfo{person}{Ulrich Trottenberg},
  \bibinfo{person}{Cornelius~W. Oosterlee}, {and} \bibinfo{person}{Anton
  Schuller}.} \bibinfo{year}{2000}\natexlab{}.
\newblock \bibinfo{booktitle}{\emph{Multigrid}}.
\newblock \bibinfo{publisher}{Elsevier}.
\newblock


\bibitem[\protect\citeauthoryear{Versteeg}{Versteeg}{1994}]%
        {versteeg1994marmousi}
\bibfield{author}{\bibinfo{person}{Roelof Versteeg}.}
  \bibinfo{year}{1994}\natexlab{}.
\newblock \showarticletitle{The Marmousi experience: Velocity model
  determination on a synthetic complex data set}.
\newblock \bibinfo{journal}{\emph{The Leading Edge}} \bibinfo{volume}{13},
  \bibinfo{number}{9} (\bibinfo{year}{1994}), \bibinfo{pages}{927--936}.
\newblock
\urldef\tempurl%
\url{https://doi.org/10.1190/1.1437051}
\showDOI{\tempurl}


\bibitem[\protect\citeauthoryear{Whigham et~al\mbox{.}}{Whigham
  et~al\mbox{.}}{1995}]%
        {whigham1995grammatically}
\bibfield{author}{\bibinfo{person}{Peter~A. Whigham} {et~al\mbox{.}}}
  \bibinfo{year}{1995}\natexlab{}.
\newblock \showarticletitle{Grammatically-based genetic programming}. In
  \bibinfo{booktitle}{\emph{Proceedings of the workshop on genetic programming:
  from theory to real-world applications}}, Vol.~\bibinfo{volume}{16}.
  \bibinfo{pages}{33--41}.
\newblock


\bibitem[\protect\citeauthoryear{Xu and Zikatanov}{Xu and Zikatanov}{2017}]%
        {xu2017algebraic}
\bibfield{author}{\bibinfo{person}{Jinchao Xu} {and} \bibinfo{person}{Ludmil
  Zikatanov}.} \bibinfo{year}{2017}\natexlab{}.
\newblock \showarticletitle{Algebraic multigrid methods}.
\newblock \bibinfo{journal}{\emph{Acta Numerica}}  \bibinfo{volume}{26}
  (\bibinfo{year}{2017}), \bibinfo{pages}{591–721}.
\newblock
\urldef\tempurl%
\url{https://doi.org/10.1017/S0962492917000083}
\showDOI{\tempurl}


\end{thebibliography}

\end{document}